\documentclass{iig}

\mathchardef\tnode="020E % temporary node
\def\arc{% unlabelled stroke for projective plane
  \hbox{\kern -0.15em
  \vbox{\hrule width 3em height 0.6ex depth -0.5 ex}
  \kern -0.33em}}

\def\darc{% double arc for GQ
  \rlap{\lower0.2ex\arc}{\raise0.2ex\arc}}

\def\stroke#1{% labelled stroke for diagrams
  \kern 0.05em
  \rlap\arc{{\textstyle{#1}}\atop\phantom\arc}
  \kern -0.22em}

\def\dstroke#1{% probably just for triple cover of Sp(4,2) quadrangle
  \kern 0.05em
  \rlap\darc{{\textstyle{#1}}\atop\phantom\darc}
  \kern -0.22em}

\def\centerscript#1{% centres text over or under a node
  \setbox0=\hbox{$\tnode$}
  \hbox to \wd0{\hss$\scriptstyle{#1}$\hss}}

\def\node{%labelled node; usage \node or \node^{label} or \node_{label}
  \def\super{}
  \def\sub{}
  \futurelet\next\dolabellednode}

  \let\sp=^
  \let\sb=_

  \def\dolabellednode{%
    \ifx\next\sb\let\next\getsub
    \else
      \ifx\next\sp\let\next\getsuper
      \else\let\next\donode
      \fi
    \fi
    \next}

  \def\getsub_#1{\def\sub{#1}\futurelet\next\dolabellednode}
  \def\getsuper^#1{\def\super{#1}\futurelet\next\dolabellednode}

  \def\donode{%
   \rlap{$\mathop{\phantom\tnode}\limits_{\centerscript{\sub}}^{\centerscript{\super}}$}\tnode}

\def\varcdn{%vertical arc with node at bottom - use as subscript to node
  \kern -0.03em\vbox{\kern -0.5ex  
  \hbox to \wd0{\hss\vrule width 0.04em depth 5.8ex\hss}
  \kern -0.3ex  \hbox{$\tnode$}}}

{

\newtheorem{phan}{Phan's Theorem}
\newtheorem{phantype}{Phan-type Theorem}
\newtheorem{solomontits}{Solomon--Tits-type Theorem}
\newtheorem{curtistits}{Curtis--Tits Theorem Version}
\newtheorem{locrec}{Local Recognition Theorem}

}

{

}
\newcommand\form{{(\cdot,\cdot)}}
\newcommand\fform{{(\form)}}

\newcommand{\mc}[1]{\mathcal{#1}}

\newcommand{\lb}{\left\{}
\newcommand{\lbr}{\left\lbrack}
\newcommand{\rb}{\right\}}
\newcommand{\rbr}{\right\rbrack}
\newcommand{\gen}[1]{\left< #1 \right>}

\newcommand{\Rad}{{\rm Rad}}
\newcommand{\Aut}{{\rm Aut}}

\newcommand{\F}{\mathbb{F}}

\newcommand{\cA}{\mathcal{A}}

\newcommand{\Fqsq}{\mathbb{F}_{q^2}}

\newcommand{\cB}{\mathcal{B}}
\newcommand\cS{{\cal S}}
\newcommand\cT{{\cal T}}
\newcommand\cC{{\cal C}}

\newcommand{\la}{\langle}
\newcommand{\ra}{\rangle}

\newcommand{\ie}{i.e., }
\newcommand{\A}{{\mc A}}

\newcommand{\G}{{\mc G}}

\newcommand{\typ}{{\rm typ}}

\newcommand\Opp{{\rm Opp}}
\newcommand\Sym{{\rm Sym}}

\newcommand\bW{\mathbf{W}}

\newcommand{\SL}{{\rm SL}}
\newcommand{\PSL}{{\rm PSL}}
\newcommand{\SU}{{\rm SU}}

\newcommand{\Sp}{{\rm Sp}}

\newcommand{\Spin}{{\rm Spin}}

\newcommand{\SO}{{\rm SO}}

\newcommand{\PSOmega}{{\rm PS}\Omega}

\newcommand\dl{\delta}

\newcommand\Sg{\Sigma}

\newcommand\eps{\epsilon}

\def\aff{\mathop{\mathrm{aff}} \nolimits}
\def\proj{\mathop{\mathrm{proj}}\nolimits}
\def\Inv{\mathop{\mathrm{Inv}}\nolimits}

\begin{document}

\title{Developments in finite Phan theory}
\author[Gramlich]{Ralf Gramlich}
\thanks{The author gratefully acknowledges his
Heisenberg fellowship of the Deutsche Forschungsgemeinschaft.}
\address{TU Darmstadt \\
FB Mathematik \\
Schlo\ss gartenstra\ss e 7 \\
64289 Darmstadt \\
Germany}
\address{The University of Birmingham \\
School of Mathematics \\
Edgbaston \\
Birmingham \\
B15 2TT \\
United Kingdom}
\email{gramlich@mathematik.tu-darmstadt.de}
\email{ralfg@maths.bham.ac.uk}
\website{http://mathematik.tu-darmstadt.de/~gramlich}

\maketitle


\section{Introduction}

Geometric methods in the theory of Chevalley groups and their generalisations have made tremendous
advances during the last few decades. Among the most noteworthy and influential of these advances are the systematic application of the concept of amalgams based on \cite{Delgado/Goldschmidt/Stellmacher:1985}, \cite{Goldschmidt:1980}, \cite{Serre:1977}, the local-to-global approach \cite{Tits:1981},
% in the Coxeter Festschrift \cite{Davis/Grunbaum/Sherk:1981}, 
ingenious applications of combinatorial topology and geometric group theory (as in \cite{Abramenko:1996}, \cite{Ivanov:1999}, \cite{Ivanov/Shpectorov:2002}, \cite{Pasini:1985}, \cite{Tits:1986}), the theory of abstract root groups \cite{Timmesfeld:1991}, \cite{Timmesfeld:1999}, \cite{Timmesfeld:2001}, and the interaction of Kac--Moody groups and twin buildings \cite{Caprace:2009}, \cite{Caprace/Muehlherr:2005}, \cite{Caprace/Muehlherr:2006}, \cite{Caprace/Remy}, \cite{Caprace/Remy:2006}, \cite{Tits:1987}, \cite{Tits:1992}.
These methods have proven fruitful over and over again in
proving, simplifying and generalising several results in group theory and have had their impact in other areas of mathematics.

The present survey attempts to give a report on the results and on the developments in recent years and to serve as a guide to the literature for the project called Curtis--Phan--Tits Theory (or, short, Phan Theory). This project has been initiated in \cite{Bennett/Gramlich/Hoffman/Shpectorov:2003} with the goal to revise Phan's results \cite{Phan:1977}, \cite{Phan:1977a} on presentations of twisted forms of finite Chevalley groups via rank one and rank two groups in order to make them accessible for the ongoing revision of the classification of the finite simple groups \cite{Gorenstein/Lyons/Solomon:1994}, \cite{Gorenstein/Lyons/Solomon:1995}, \cite{Gorenstein/Lyons/Solomon:1998}, \cite{Gorenstein/Lyons/Solomon:1999}, \cite{Gorenstein/Lyons/Solomon:2002}, \cite{Gorenstein/Lyons/Solomon:2006}. 

The main impact of Phan's results \cite{Phan:1977}, \cite{Phan:1977a} in the classification can be seen in \cite{Aschbacher:1977} side by side with the famous Curtis--Tits Theorem established in \cite{Curtis:1965}, \cite{Tits:1962}, \cite[Theorem 13.32]{Tits:1974}; see also \cite[Section 2.9]{Gorenstein/Lyons/Solomon:1998} and Section \ref{Curtis--Tits theorem} of this survey. 
As in this survey my main concern is the revision of Phan's results, I refer to \cite{Aschbacher:2004} regarding the current overall state of the classification of the finite simple groups. 

Phan Theory enters the stage in what is called Step 2 in \cite{Aschbacher:2004}, the identification of the {\em minimal
counterexample} $G$ as one of the known simple groups.  By Step 1, the local analysis, inside the minimal
counterexample $G$ one reconstructs one or more of the proper
subgroups using the inductive assumption and available techniques.
Thus, the initial point of the identification, Step 2,  is a set of subgroups of
$G$ that resemble the subgroups of a central extension $\widehat{G}$ of
some known simple group, referred to as the {\em target group};  the
output of the identification step is the statement that $G$ is
isomorphic to a central quotient of $\widehat{G}$.  
Two of the most widely used
identification tools in this step are the Curtis--Tits Theorem and Phan's theorems. I have already mentioned \cite{Aschbacher:1977} as one of the main occurrences of these tools in the classification of the finite simple groups and refer the reader to \cite[Section 7.5]{Gorenstein/Lyons/Solomon:2006} for an occurrence of Phan's revised results in the revision of the classification. In Section \ref{app} I describe a possible setup for an application of Phan Theory in the revision of the classification via centralisers of involutions.

The Curtis--Tits Theorem allows the identification of $G$ with a quotient of a universal
Chevalley group $\widehat{G}$ of twisted or untwisted type provided that
$G$ contains a generating system of subgroups identical to the system of
appropriately chosen rank two Levi factors of $\widehat{G}$. For instance,
in the case of the Dynkin diagram $A_n$, the system in question consists of
all the groups $\SL_2(\F)$, $\SL_3(\F)$, and $\SL_2(\F)\times \SL_2(\F)$ lying block-diagonally in
$\widehat{G} \cong \SL_{n+1}(\F)$, considered as a matrix group with respect to a suitable basis of its natural module. Phan's first theorem \cite{Phan:1977} on the other hand deals with the case $\widehat{G}\cong\SU_{n+1}(q^2)$, considered as a matrix group with respect to an orthonormal basis of its natural module, and the
system of block-diagonal subgroups $\SU_2(q^2)$, $\SU_3(q^2)$, and  $\SU_2(q^2)\times
\SU_2(q^2)$. 

In this survey I will describe how a systematic geometric approach making serious use of buildings and twin buildings yields a Phan-type theorem for Chevalley groups of each irreducible spherical type of rank at least three. The complete result is stated in Section \ref{complete}.

\subsection*{Acknowledgements}

The author thanks the main protagonists of Phan theory, Curt Bennett, Corneliu Hoffman, Bernhard M\"uhlherr, and Sergey Shpectorov, without whom this survey would have nothing to report. The author is particularly grateful for the fact that he has been invited into the project by Curt Bennett, Corneliu Hoffman, and Sergey Shpectorov as a visiting PhD student at Bowling Green State University in 2001. Moreover, the author expresses his deep gratitude to Richard Lyons and Ronald Solomon without whose question for a revision of Phan's results there would not have been any research in this direction. Furthermore, the author very gratefully acknowledges Michael Aschbacher's, Antonio Pasini's, and the referee's extensive
comments on earlier versions of this article. Additional very helpful remarks and suggestions have been made by Pierre-Emmanuel Caprace, Alice Devillers, and Franz Timmesfeld. Finally, the author expresses his gratitude to Max Horn for a careful proofreading.

The author and part of his research group are supported by the Deutsche Forschungsgemeinschaft via the grants GR 2077/4, GR 2077/5, and GR 2077/7.

{\footnotesize
\addcontentsline{toc}{section}{Contents}
\setcounter{tocdepth}{2}
\tableofcontents
}

\section{Geometries and amalgams} \label{geometries}

In this section I give a quick overview over the basic geometric notions and results used in the present survey. Most of these notions are due to Tits \cite{Tits:1981}.
I refer the reader to \cite{Tits:1981} and also to \cite{Buekenhout:1995}, \cite{Buekenhout/Cohen}, \cite{Pasini:1994} for a well-founded introduction to synthetic geometry including proofs and many helpful examples.

\subsection{Geometries}

\subsubsection{Pregeometries and geometries}

A {\it pregeometry} $\G = (X,*,\typ)$ over a finite set $I$ is a set of elements $X$ together with a {\em type function} $\typ : X \rightarrow I$ and a reflexive and symmetric {\em incidence relation}
$* \subseteq X \times X$  
such that for
any two elements $x, y \in \G$ with $x * y$ and $\typ(x) = \typ(y)$ we
have $x=y$.  The {\it rank} of a pregeometry $\G$ is the cardinality of its type set $I$. A {\it flag} in $\G$ is a set of pairwise incident
elements. Hence the type function injects any flag into the type
set, this image is called the {\em type} of the flag. A {\it geometry} is a pregeometry with the property that $\typ$ induces a
bijection between any maximal flag of $\G$ and $I$. Flags of type $I$ are called {\em chambers}.
\subsubsection{Residues}

The {\it residue} $\G_F$ of a flag $F$ in a pregeometry $\G$ consists of
the set of elements from $\G \setminus F$ that are incident to all
elements of $F$ with the restricted incidence and type functions, the latter co-restricted to $I \backslash \typ(F)$, turning the residue $\G_F$ into a
pregeometry over $I \setminus \typ(F)$. If $\G$ is a geometry over $I$, then any of its residues $\G_F$ is a geometry over $I \backslash \typ(F)$. 
  The rank of the residue of a
flag $F$ is called the {\it co-rank} of $F$.  
A non-empty pregeometry $\G$ is
said to be {\it connected}, if the graph $(X,*)$ is connected.  
Following \cite[Section 1.2]{Tits:1981}, a pregeometry $\G$ is {\it residually
connected}, if the residue in $\G$ of any flag of co-rank at least two
is connected and the residue of any flag of
co-rank one is non-empty. A residually connected pregeometry is automatically a geometry.

\subsubsection{Automorphisms}

An {\em automorphism}\index{automorphism} of a pregeometry $\G$ over $I$ is a permutation of its elements
that preserves type and incidence and whose inverse permutation also preserves incidence.  The group of all automorphisms of
$\G$ will be denoted by $\Aut\ \G$.  A subgroup $G \le \Aut\ \G$ acts {\it flag-transitively} on $\G$ if, for each $J \subseteq I$, it is transitive on the set of all
flags of type $J$. A group $G$ of automorphisms of a geometry $\G$ over $I$ is flag-transitive if and only if $G$ is transitive on the set of maximal flags of $\G$, because each flag of $\G$ can be extended to a flag of type $I$ of $\G$. A pregeometry that admits a flag-transitive automorphism
group is called flag-transitive.
A {\em parabolic subgroup} (or simply a {\em parabolic}) of $G$ is the
stabiliser in $G$ of a non-empty flag $F$ of $\G$.  The {\em rank}\index{rank} of
the parabolic is defined as the co-rank of $F$.

The term parabolic subgroup is inspired by the parabolic subgroups of algebraic groups which occur as stabilisers of residues of buildings.

\subsection{Simplicial complexes} \label{simplicial}

There exist very many good books dealing with the theory of simplicial complexes, many of them with very different flavours, ranging from combinatorics and graph theory \cite{Gross/Tucker:1987} to differential geometry \cite{Bredon:1993} and to topology \cite{Spanier:1966}. Other classical references are \cite{Seifert/Threlfall:1934}, \cite{tomDieck:2000}.

\subsubsection{Complexes}

A {\em simplicial complex}\index{simplicial complex}\index{complex!simplicial} $\cS$ is a pair $(X,\Delta)$ where $X$ is a
set and $\Delta$ is a collection of non-empty finite subsets of $X$ containing each subset of $X$ of cardinality one such that
$A\in\Delta$ and $\emptyset \neq B\subseteq A$ implies $B\in\Delta$.  The elements of
$\Delta$ are called {\em simplices}. A simplicial complex in which each chain of simplices is finite is called {\it pure}, if all of its maximal simplices have the same cardinality.

A {\em morphism}\index{morphism} from a complex $\cS=(X,\Delta)$ to a complex
$\cS'=(X',\Delta')$ is a map between $X$ and $X'$ which takes
simplices to simplices.  The {\em star}\index{star} of a simplex $A\in\Delta$ is
the set of subsets $B\in\Delta$ such that $A\subseteq B$.  A {\em
covering}\index{covering} is a surjective morphism $\phi$ from $\cS$ to $\cS'$ such that for
every $A\in\Delta$, the function $\phi$ maps the star of $A$
bijectively onto the star of $\phi(A)$.
A {\em path} on a complex $\cS$ is a finite sequence $x_0$, $x_1$, \ldots, $x_n$
of elements of $X$ such that $x_{i-1}$ and $x_i$ are contained in a
simplex for all $i=1, \ldots, n$. The complex $\cS$ is {\em connected}, if
every two elements of $X$ can be connected by a path.  

\subsubsection{Homotopy} \label{222}

The following three operations are called {\em elementary homotopies} of paths:  substituting a subsequence $x$, $x$ (a {\em repetition}\index{repetition}) by $x$, or vice versa,  substituting a subsequence
$x$, $y$, $x$ (a {\em return}\index{return}) by $x$, or vice versa, or  substituting a subsequence
$x$, $y$, $z$, $x$ (a {\em triangle}\index{triangle}) by $x$ or vice versa, provided that $x$, $y$, $z$ form a
simplex. 
 Two paths are {\em homotopically equivalent} if they can be obtained
from one another in a finite sequence of elementary homotopies.  A
{\em cycle}, that is, a path with $x_0=x_n$, is called {\em null-homotopic}\index{path!null-homotopic}, if it is homotopically equivalent to the trivial path
$x_0$.  The {\em fundamental group}\index{fundamental group}\index{group!fundamental} $\pi_1(\cS, x)$, where $x \in X$, is
the set of homotopy classes of cycles based at $x$ where the product is defined to be concatenation of cycles.
The fundamental group is independent of the choice of
the base vertex $x$ inside a fixed connected component, while it may vary for base vertices in distinct
connected components. When considering connected complexes only, the coverings of $\cS$, taken up to a certain
natural equivalence, correspond bijectively to the subgroups of
$\pi_1(\cS, x)$, cf.\ \cite[\S 55]{Seifert/Threlfall:1934}.  A connected complex $\cS$ is called {\em simply
connected}, if each covering $\cS' \to \cS$ with connected $\cS'$ is an isomorphism, or equivalently (\cite[\S 56]{Seifert/Threlfall:1934}, \cite[Section 2.5]{Spanier:1966}), if
$\pi_1(\cS, x)=1$.

\subsubsection{Flag complexes and realisations} \label{223}

With every pregeometry $\G = (X,\typ,*)$ one can associate its {\em flag complex} which is a simplicial complex defined on the set $X$ whose simplices
are the flags of $\G$. The flag complex of a pregeometry $\G$ is pure if and only if $\G$ is a geometry. A pregeometry $\G$ is {\em simply
connected}, if such is its flag complex.

For a simplicial complex $\cS=(X,\Delta)$ denote by $|\cS|$ the set of all functions $\alpha$ from $X$ to the real unit interval $I$ satisfying that the set $\{ v \in X \mid \alpha(v) \neq 0 \}$ is contained in $\Delta$ and that $\sum_{v \in X} \alpha(v) = 1$, i.e., $|\cS|$ is obtained from $\cS$ via barycentric coordinates. In this survey I consider the weak (coherent) topology on $|\cS|$, cf.\ \cite[3.1.14]{Spanier:1966}, and call it the {\em realisation} of $\cS$. With respect to this topology, the fundamental group $\pi_1(\cS,x)$ defined combinatorially in Section \ref{222} coincides with the usual fundamental group defined topologically, see \cite[3.6.17]{Spanier:1966}. 

%Let $\cS=(X,\Delta)$ be a simplicial complex and let $s$ be a simplex. The {\em closed simplex} $|s|$ equals $\{ \alpha \in |\cS| \mid \alpha(v) \neq 0 \Longrightarrow v \in s \}$. Hence, if $s$ is a $q$-simplex, then $|s|$ is in one-to-one correspondence with the set $\{ x \in \mathbb{R}^{q+1} \mid 0 \leq x_i \leq 1, \sum x_i = 1 \}$. The {\em open simplex} $\gen{s}$ equals $\{ \alpha \in |\cS| \mid \alpha(v) \neq 0 \Longleftrightarrow v \in s \}$. The {\em open simplex} $\gen{s}$ is open in $|s|$, but not necessarily in $|\cS|$.
%For a vertex $v$ of $\cS$, the {\em open star} $St_\cS(v)$ equals $\{ \alpha \in |\cS| \mid \alpha(v) \neq 0 \}$ and for a simplex $s$ the {\em open star} $\mathrm{St}_\cS(s)$ equals $\{ \alpha \in |\cS| \mid \alpha(v) \neq 0 \Longleftarrow v \in s \}$. If $s$ is a maximal simplex, then $\mathrm{St}_\cS(s) = \gen{s}$. More generally, $\mathrm{St}_\cS(s) = \bigcup \{ \gen{t} \mid s \subseteq t \}$, cf.\ \cite[3.1.24]{Spanier:1966}.
%The {\em link} of $s$ is defined as $\mathrm{Lk}_\cS(s) := \overline{\mathrm{St}_\cS(s)} \backslash \mathrm{St}_\cS(s)$.

\subsubsection{Wedges of spheres}

Let $X$ and $Y$ be pointed spaces, i.e., topological spaces with distinguished base points $x_0$ and $y_0$. Then the {\em wedge sum} $X \vee Y$ of $X$ and $Y$ is the quotient of the disjoint union $X \sqcup Y$ by the identification $x_0 \sim y_0$, i.e., $$X \vee Y := (X \sqcup Y) / \{ x_0 \sim y_0 \}.$$

In general, if $(X_i)_{i \in I}$ is a family of pointed spaces with base points $(x_i)_{i \in I}$, then the {\em wedge sum} of this family is given by $$\bigvee_{i \in I} X_i := \bigsqcup_{i \in I} X_i / \{ x_i \sim x_j \mid i, j \in I \}.$$ 
The wedge sum of a family of spheres of the same dimension $n$ is called a {\em wedge of spheres} or, if one wants to specify the dimension, a {\em wedge of $n$-spheres}.

\subsection{Chamber systems} \label{chsys}

Chamber systems and their interaction with pregeometries and simplicial complexes as introduced in \cite{Tits:1981} play a crucial role in this survey. Details on chamber systems can also be found in \cite{Abramenko/Brown:2008}, \cite{Brown:1989}, \cite{Buekenhout/Cohen}, \cite{Pasini:1994}, \cite{Ronan:1989}, \cite{Weiss:2003}. In this section I sketch the most fundamental information and try to highlight some interaction with the objects introduced before.

\subsubsection{Chambers}\label{2.3.1}

A {\it chamber system} $\mathcal{C} = (C,(\sim_i)_{i \in I})$ over a type set $I$ is a set $C$, called the
set of {\em chambers},\index{chamber} together with equivalence relations $\sim_i$, $i \in
I$.  For $i \in I$ and chambers $c, d \in \cC$, the chambers $c$ and
$d$ are called {\em $i$-adjacent} if $c\sim_i d$. The chambers $c$, $d$ are {\em adjacent} if they are $i$-adjacent for some $i \in I$.

A chamber system $\cC$ is called {\it thick} if for every $i\in I$
and every chamber $c\in \cC$, there are at least three chambers ($c$
and two other chambers) $i$-adjacent to $c$; it is called {\em thin} if for every $i\in I$
and every chamber $c\in \cC$, there are exactly two chambers ($c$
and one other chamber) $i$-adjacent to $c$.  

A {\it gallery} in ${\cal C}$ is a finite sequence
$(c_0,c_1,\ldots,c_t)$  such that $c_{k} \in C$ 
for all $0 \leq k \leq t$ and such that
$c_{k-1}$ is adjacent to $c_{k}$ for all
$1 \leq k \leq t$. The number $t$ is called the {\it length} of the
gallery. 
The chamber system ${\cal C}$ is said to 
be {\it connected}, if for any two chambers there exists
a gallery joining them.

For $J \subseteq I$, the {\em $J$-residue} of a chamber $c$ is the chamber system $\mathcal{R}_J(c) = (R_J(c),(\sim_j)_{j \in J})$ consisting of those chambers of $\cC$ that can be connected to $c$ via a gallery using $j$-adjacencies ($j \in J$) only; such galleries are called {\it $J$-galleries}. A $J$-residue with $|J| = 1$ is called a {\em panel}.

\subsubsection{Chamber systems and pregeometries} \label{csap}

If $\G$ is a pregeometry with type set $I$, then one can construct a
chamber system $\cC=\cC(\G)$ over $I$ as follows.  The chambers are
the flags of $\G$ of type $I$ and two such flags are $i$-adjacent if and
only if they contain the same element of type $j$ for all $j \in
I\setminus\{i\}$.  A chamber system is called {\it geometric}\index{chamber system!geometric}, if it
can be obtained in this way.

If, conversely, $\mathcal{C}$ is a chamber system over $I$, the {\em pregeometry of $\mathcal{C}$} (denoted by $\G(\mathcal{C})$) is the pregeometry over $I$ whose elements of type $i$ are the pairs $(x,i)$ with $x$ an $I \backslash \{ i \}$-residue of $\mathcal{C}$ in which two elements $(x,k)$, $(y,l)$ of $\G(\mathcal{C})$ are incident if and only if $x \cap y \neq \emptyset$ in $\mathcal{C}$, cf.\ \cite{Tits:1981}.
If $\psi_\mathcal{C}(c)$, for $c \in \mathcal{C}$, denotes the set of all $I \backslash \{ i \}$-residues, $i \in I$, containing $c$, then the map $$\mathcal{C} \to \mathcal{C}(\G(\mathcal{C})) : c \mapsto \psi_\mathcal{C}(c)$$ is a homomorphism of chamber systems, by \cite[Proposition 3.5.6]{Buekenhout/Cohen}.

In general, $\G \not\cong \G(\cC(\G))$ and $\cC \not\cong \cC(\G(\cC))$, see \cite[Section 2.2]{Tits:1981}. However, if $\G$ is residually connected, then $\G \cong \G(\cC(\G))$, cf.\ \cite[Section 3.5]{Buekenhout/Cohen}, \cite[Section 2.2]{Tits:1981}. Moreover,  by \cite[Theorem 3.5.7]{Buekenhout/Cohen}, the homomorphism $$\mathcal{C} \to \mathcal{C}(\G(\mathcal{C})) : c \mapsto \psi_\mathcal{C}(c)$$ is an isomorphism if and only if
for any set $\{ (x_i,i) \mid i \in I \}$ (where $x_i$ is an $I \backslash \{ i \}$-residue of $\mathcal{C}$) such that $x_i \cap x_j \neq \emptyset$ for all $i, j \in I$, the intersection $\bigcap_{i \in I} x_i$ is non-empty, and for distinct chambers $c$, $d$ of $\mathcal{C}$ there is some $I \backslash \{ i \}$-residue of $\mathcal{C}$ containing $c$ but not $d$.

\subsubsection{Homotopy}\label{2-connected}

The concept of homotopy introduced for simplicial complexes, cf.\ Section \ref{simplicial}, can also be defined for chamber systems. Excellent sources are \cite{Pasini:1994}, \cite{Tits:1981}.

Let
 $m \geq 1$ be an integer and let
$({\cal C},(\sim_i)_{i \in I})$ be a chamber system over a set $I$.
Two galleries $G =(c_0,\ldots,c_k)$ and 
$H = (c_0',\ldots,c_{k'}') $ are called
{\it elementarily $m$-homotopic}, if
there exist two galleries $X,Y$ 
and two $J$-galleries $G_0,H_0$ for some $J \subset I$
of cardinality at most $m$ 
such that $G = XG_0Y$, $H=XH_0Y$. 
Two galleries $G$, $H$ are said to be {\it $m$-homotopic} if
there exists a finite sequence $G_0$, $G_1$, \ldots, $G_l$ of  galleries 
such that $G_0 = G$, $G_l = H$ and such that
$G_{k-1}$ is elementarily $m$-homotopic to $G_{k}$ for all
$1 \leq k \leq l$.
A closed gallery $G$ is called {\it null-$m$-homotopic}
if it is $m$-homotopic to the gallery consisting of the initial chamber of $G$. 

The chamber system ${\cal C}$ is called
{\it simply $m$-connected}, if it is connected and if each closed
gallery is null-$m$-homotopic.
Given a gallery $G$, then $GG^{-1}$ is null-$m$-homotopic.
Furthermore, two galleries $H$, $G$ are $m$-homotopic if and only if
the gallery $GH^{-1}$ is null-$m$-homotopic.

If $\mathcal{C}$ is a chamber system over a finite set $I$ such that the map $\mathcal{C} \to \mathcal{C}(\G(\mathcal{C})) : c \mapsto \psi_\mathcal{C}(c)$ from Section \ref{csap} is an isomorphism, then $\mathcal{C}$ is simply $(|I|-1)$-connected if and only if $\G(\mathcal{C})$ is simply connected. For $m < |I| -1$ it is unknown to me what it means for a geometry, if its chamber system is simply $m$-connected. Note that there exists a rank four geometry (cf.\ \cite{Stroth:1989}) for McLaughlin's sporadic simple group $McL$ whose chamber system is simply $2$-connected and which admits residues of rank three which are not simply connected.

\subsection{Coset pregeometries and reconstruction}

\subsubsection{Coset pregeometries}

Let $I$ be a set, let $G$ be a group and let $(G_{i})_{i \in I}$ be a family of subgroups of $G$. Then $$(\bigsqcup_{i \in I} G/G_{i},*,\typ)$$ with $\typ (gG_{i}) = i$ and $gG_{i} * hG_{j}$ if and only if $gG_{i} \cap hG_{j} \neq \emptyset$ is a pregeometry of type $I$, the {\em coset pregeometry of $G$} with respect to $(G_{i})_{i \in I}$. The groups $G_{i}$ are called the {\em maximal parabolic subgroups} of the coset pregeometry. Since the type function is completely determined by the indices, we denote the coset pregeometry of $G$ with respect to $(G_{i})_{i \in I}$ by $((G/G_{i})_{i \in I}, *)$. 
The family $(G_i)_{i \in I}$ forms a chamber of the coset pregeometry, called the {\em base chamber}. For $J\subseteq I$ define
$G_J:=\bigcap_{j\in J}G_j$.

\subsubsection{Reconstruction} \label{trans}

Certainly any coset pregeometry is {\em incidence-transitive}, i.e., for any two flags $c$ and $d$ with $|\typ(c)| = 2 = |\typ(d)|$ and $\typ(c) = \typ(d)$ there exists an element $g \in G$ that maps $c$ onto $d$. Indeed, if $gG_{i} \cap hG_{j} \neq
\emptyset$, then choose $a \in gG_i \cap hG_j$. It follows $aG_i = gG_i$ and $aG_j = hG_j$ and therefore the automorphism $a^{-1}$
maps the incident pair $gG_{i}$, $hG_{j}$ onto the incident pair
$G_{i}$, $G_{j}$.
Conversely, any incidence-transitive pregeometry can be described as a coset pregeometry via its parabolic subgroups. 

If $\G = (X,*,\typ)$ is a pregeometry over $I$ with an incidence-transitive group $G$ of automorphisms of $\G$ and a maximal flag $F = (x_i)_{i \in I}$ of $\G$, then the bijection $$((G/G_{x_i})_{i \in I},*') \rightarrow \G : g G_{x_i} \mapsto gx_i$$ is an isomorphism between pregeometries and between $G$-sets. (Recall here that two actions $\phi
: G \to \Sym\ M$ and $\phi' : G \to \Sym\ M'$ are called
{\em isomorphic},
if there is a bijection
$\psi : M \to M'$ such that
$\psi \circ \phi(g) \circ \psi^{-1} = \phi'(g)$ for each $g\in G$ or,  equivalently,
$\psi \circ \phi(g) = \phi'(g) \circ \psi$ for all $g\in G$; in this case,
we also say that $M$ and $M'$ are {\em isomorphic $G$-sets}.) The observation of this isomorphism $((G/G_{x_i})_{i \in I},*') \rightarrow \G : g G_{x_i} \mapsto gx_i$ goes back to \cite{Klein:1872} and has been proved formally in \cite{Freudenthal:1951}, \cite{Freudenthal:1985}.

It happens quite frequently that interesting geometries are not incidence-transitive. This is also the case in Phan theory, see e.g.\ Section \ref{flippol}, so that often a more general definition of a coset pregeometry is necessary. I refer the reader to \cite{Gramlich/Horn/Pasini/Maldeghem}, \cite{Gramlich/Maldeghem:2006}, \cite{Stroppel:1992}, \cite{Stroppel:1993}, \cite{Wich:1996} for details. The most general concept of reconstruction in this context known to me are complexes of groups as treated in \cite[Chapter III.$\mathcal{C}$]{Bridson/Haefliger:1999}.

\subsection{Geometric covering theory and Tits' Lemma} \label{titslemmasec}

\subsubsection{Amalgams}

An {\it amalgam} $\mathcal{A}$ of groups is a set with a
partial operation of multiplication and a collection of subsets $\{ G_i\}_{i \in I}$, for some index set $I$, such that (i) $\mathcal{A} = \bigcup_{i \in I} G_i$, (ii) for each $i \in I$, the restriction of the
multiplication to $G_i$ turns $G_i$ into a group, (iii) the
product $ab$ is defined if and only if $a,b \in G_i$ for some $i\in
I$, and (iv) $G_i\cap G_j$ is a subgroup of $G_i$ and $G_j$ for all
$i,j \in I$.

An {\em enveloping group} of an amalgam $\cA$ is a group $G$ together with a
mapping $\phi$ from $\cA$ to $G$ such that the restriction of $\phi$
to every $G_i$ is a homomorphism and $\phi(\cA)$ generates $G$.  The
{\em universal enveloping group} of $\cA$ is isomorphic to the group $U(\cA)$ with generators
$\{ t_s \mid s \in \cA\}$ and relations $t_xt_y=t_{xy}$ whenever $x,y \in
G_i$ for some $i$; the corresponding mapping is given by $x\mapsto
t_x$, see \cite[Chapter I, Section 1.1, Proposition 1]{Serre:2003}. Every enveloping group is isomorphic to a quotient of
the universal enveloping group $U(\cA)$.
For more details we refer the reader to \cite{Serre:1977}, \cite{Serre:2003}. Intransitive geometries may lead to {\em fused amalgams} as defined and studied in \cite{Gramlich/Horn/Pasini/Maldeghem}. Again, I also refer to the concept of complexes of groups \cite[Chapter III.$\mathcal{C}$]{Bridson/Haefliger:1999}.

In some references an enveloping group of an amalgam is called a completion of this amalgam.

\subsubsection{The fundamental theorem of geometric covering theory} \label{ft}

Many identification problems in group theory amount to finding the universal enveloping groups of certain
amalgams arising inside some abstract group, for instance as stabilisers of some group action on some simplicial complex with a fundamental domain. The result that connects such
amalgams and their enveloping groups with combinatorial-topological properties of the set acted on is a lemma proved in \cite{Pasini:1985}, \cite{Tits:1986}, cf.\ Section \ref{tl} known as Tits' lemma. It can be obtained as a corollary of the Fundamental Theorem of geometric covering theory discussed in this section.

Suppose $\G$ is a geometry and $G\le\Aut\ \G$ is an incidence-transitive 
group of automorphisms.  Corresponding to $\G$ and $G$ and some maximal flag $F$, there is an amalgam 
$\A=\A(\G,G,F)$, the {\em amalgam of parabolics with respect to $\G$, $G$, $F$},\index{amalgam!of parabolics} defined as the family $(G_E)_{\emptyset \neq E \subseteq F}$, where $G_E$ 
denotes the stabiliser of $\emptyset \neq E \subseteq F$ in $G$.  

In case $G$ is flag-transitive, 
the amalgam $\A$ is independent up to conjugation of the choice 
of $F$. If $\G$ is connected, then $\cA$ generates $G$ (\cite[Lemma 1.4.2]{Ivanov:1999}), so that
$G$ is an enveloping group of $\cA$.
One of the main tools for geometric proofs of group-theoretic identification theorems is the Fundamental Theorem of geometric covering theory, see \cite{Ivanov/Shpectorov:2002}.  

\medskip
\noindent {\bf Fundamental Theorem of geometric covering theory (\cite[Theorem 1.4.5]{Ivanov/Shpectorov:2002}).}
{\em Let $\mc{G} = (X,*,\typ)$ be a connected geometry over $I$ of rank at least three, and let $G$ be a flag-transitive group of automorphisms of $\G$. Moreover, let $F$ be a maximal flag and let $\mc{A} = \mc{A}(\G,G,F)$ be the corresponding amalgam of parabolics. Then the coset pregeometry $\widehat\G = ((\mc{U}(\mc{A})/G_{x})_{x \in F},*)$ is a simply connected geometry that admits a covering $\pi : \widehat\G \rightarrow \G$ induced by the natural epimorphism $\mc{U}(\mc{A}) \rightarrow G$. Moreover, $\mc{U}(\mc{A})$ is of the form $\pi_1(\G).G$, i.e., $\mathcal{U}(\mc{A})/\pi_1(\G) \cong G$.}

\subsubsection{Tits' Lemma} \label{tl}

An immediate consequence of the Fundamental Theorem is Tits' Lemma, cf.\  \cite[Corollary 1.4.6]{Ivanov/Shpectorov:2002}, \cite[Lemma 5]{Pasini:1985}, \cite[Theorem 12.28]{Pasini:1994}, \cite[Corollary 1]{Tits:1986}.

\medskip
\noindent {\bf Tits' Lemma (\cite[Corollary 1.4.6]{Ivanov/Shpectorov:2002}).} {\em Let $\mc{G} = (X,*,\typ)$ be a connected geometry over $I$ with a flag-transitive group $G$ of automorphisms of $\G$, let $F$ be a maximal flag $\G$, and let $\mc{A}(\G,G,F)$ be the corresponding amalgam of parabolics. Then the geometry $\G$ is simply connected if and only if the canonical epimorphism $\mc{U}(\mc{A}(\G,G,W)) \rightarrow G$ is an isomorphism.}

\medskip
This result reduces the problem of identifying the universal
enveloping group of a certain amalgam to proving that the corresponding
geometry is simply connected, i.e., proving that the fundamental group of
its flag complex is trivial. 
Geometric covering theory has been extended to certain classes of intransitive geometries, leading to more general concepts of amalgams and different versions of the Fundamental Theorem and Tits' Lemma. I refer the reader to \cite{Gramlich/Horn/Pasini/Maldeghem} and \cite{Gramlich/Maldeghem:2006} for details. Again also the concept of complexes of groups \cite[Chapter III.$\mathcal{C}$]{Bridson/Haefliger:1999} should be mentioned.

\subsubsection{Shapes} \label{redundancy}

For a maximal flag $F$ a {\em shape} is a subset $\mc{W}$ of $2^F$ such that $2^F \ni U' \supset U \in \mc{W}$ implies $U' \in \mc{W}$, i.e., $\mc{W}$ is a subset of the power set of $F$ that is closed under passing to supersets.
The {\em amalgam of shape $\mc{W}$ with respect to $\mc{G}$, $G$, $F$} is the family $(G_U)_{U \in \mc{W}}$, where $G_U$ is the stabiliser of $U \in \mathcal{W}$ in $G$. It is denoted by $\mc{A}_\mc{W}(\mc{G},G,F)$. 
Shapes allow for a neat explanation why many presentations of groups based on amalgams of parabolics are redundant. 

\medskip
\noindent {\bf Redundancy Theorem (\cite[Theorem 3.3]{Gramlich/Maldeghem:2006}).} {\em Let $\mc{G} = (X,*,\typ)$ be a geometry over some finite set $I$, let $G$ be a flag-transitive group of automorphisms of
$\G$, and let $F$ be a maximal flag of $\G$. Moreover, let $\mc{W} \subseteq 2^F$ be a shape, assume that for each flag $U \in 2^F \backslash \mc{W}$ the residue
$\G_U$ is simply connected, and let $\mc{A}(\G,G,F)$ and $\mc{A}_{\mc{W}}(\G,G,F)$ be the amalgam of maximal
parabolics, resp.\ the amalgam of shape $\mc{W}$ of $\G$ with respect to $G$ and $F$. Then $G =
\mc{U}(\mc{A}_{\mc{W}}(\G,G,F))$ and, if $\emptyset \not\in \mc{W}$, furthermore $G =
\mc{U}(\mc{A}(\G,G,F)) = \mc{U}(\mc{A}_{\mc{W}}(\G,G,F))$.
}

\section{Phan's Theorems} \label{Phan theorem}

\subsection{Phan's first theorem} \label{phancon}

The first of the group-theoretic identification theorems I discuss in this survey is Phan's
first theorem.  In 1977 Kok-Wee Phan \cite{Phan:1977} | the namesake of the theory reported on in this survey | described a method of identification of a group $G$ as a quotient of the unitary group $\SU_{n+1}(q^2)$ via a generating configuration consisting of subgroups $\SU_2(q^2)$ and $\SU_3(q^2)$ and $\SU_2(q^2)\times \SU_2(q^2)$ in $G$. 

\subsubsection{Phan systems}

It is helpful to begin by looking at this
configuration of subgroups inside $\SU_{n+1}(q^2)$ in order to motivate the forthcoming
definitions.
For $n\ge 2$ and $q$ a prime power, consider
$G=\SU_{n+1}(q^2)$ acting as a matrix group with respect to an orthonormal basis on a unitary $(n+1)$-dimensional vector space over $\Fqsq$, and let $U_i\cong \SU_2(q^2)$, $i=1,2,\ldots,n$, be the subgroups of $G$ corresponding to the $(2\times 2)$-blocks along the
main diagonal represented as matrix groups with respect to the chosen orthonormal basis.  Let $T_i$ be the diagonal subgroup in $U_i$ with respect to this basis, which is a maximal torus of $U_i$ of size $q+1$.  For $q\geq 3$ and $1 \le i, j \le n$ the subgroups $U_i$ and $T_i$ satisfy the following axioms:
\begin{description}
\item[{\bf (P1)}] if $|i-j|>1$, then $[x,y]=1$ for all $x\in U_i$ and
$y\in U_j$,
\item[{\bf (P2)}] if $|i-j|=1$, then $\la U_i,U_j\ra$ is isomorphic to
$\SU_3(q^2)$; moreover $[x,y]=1$ for all $x\in T_i$ and $y\in T_j$, and
\item[{\bf (P3)}] $G = \gen{U_i \mid 1 \leq i \leq n}$.
\end{description}

\subsubsection{Phan's Theorem} \label{phanasch}

If $G$ is an arbitrary group containing a system of subgroups
$U_i\cong \SU_2(q^2)$ with a particular maximal torus $T_i$ of size $q+1$ chosen in each $U_i$ such that the conditions (P1), (P2), (P3) hold
for $G$, then one says that $G$ admits a {\em Phan system of type $A_n$ over $\Fqsq$}. In \cite{Aschbacher:1977} this configuration is called a {\em generating system} of type $I$, the groups $U_i$ are called {\em fundamental subgroups}. In that paper the following theorem, Phan's Theorem, is applied to obtain a characterisation of Chevalley groups over finite fields of odd order; note the additions made in \cite{Aschbacher} and \cite{Aschbacher:1980}.

\begin{phan}[Phan \cite{Phan:1977}] \label{Phan}
Let $q \geq 5$, let $n \geq 3$, and let $G$ be a group admitting a Phan system of type $A_n$ over $\Fqsq$.
Then $G$ is isomorphic to a quotient of $\SU_{n+1}(q^2)$.
\end{phan}

My favourite way of proving a result like Phan's Theorem \ref{Phan} is to translate the statement into an amalgamation problem. This means that one first constructs an abstract amalgam from the Phan system and proves that up to central extensions and isomorphisms any such amalgam is unique. Second one proves that the group admitting the Phan system is a central quotient of the universal enveloping group of the constructed unique amalgam. The first step has been well understood by now, cf.\ \cite{Bennett/Shpectorov:2004}, \cite{Dunlap:2005}, also \cite{Muehlherr:1999}. Therefore in this survey I will only concern myself with the second step. 

\subsection{Aschbacher's geometry and its simple connectedness} \label{aschgeom}

\subsubsection{Weak Phan systems of type $A_n$} \label{weakphan}

I will describe the proof of a slightly more general statement than Phan's Theorem \ref{Phan}. Following \cite{Bennett/Shpectorov:2004}, a group $G$ admits a {\em weak Phan system of type $A_n$ over $\mathbb{F}_{q^2}$}, if
$G$ contains subgroups $U_i\cong \SU_2(q^2)$, $i=1,2,\ldots,n$, and
$U_{i,j}$, $1\le i<j\le n$, so that the following hold:
\begin{description}
\item[{\bf (wP1)}] if $|i-j|>1$, then $[x,y]=1$ for all $x\in U_i$ and
$y\in U_j$,
\item[{\bf (wP2)}] if $|i-j|=1$, the groups $U_i$ and $U_{j}$ are contained
in $U_{i,j}$, which is isomorphic to a central quotient of $\SU_3(q^2)$;
moreover, $U_i$ and $U_{j}$ form a standard pair (see below) in $U_{i,j}$, and
\item[{\bf (wP3)}] $G = \gen{U_{i,j} \mid 1 \leq i < j \leq n}$.
\end{description}
Here a {\it standard pair}\index{standard pair} in the matrix group $\SU_3(q^2)$ is a pair of subgroups isomorphic to
$\SU_2(q^2)$ conjugate as a pair to the two block-diagonal groups isomorphic to
$\SU_2(q^2)$, i.e., these two groups centralise a pair of orthonormal vectors of the natural module of $\SU_3(q^2)$.  Standard pairs in central quotients of $\SU_3(q^2)$ are defined as the
images under the canonical homomorphism of standard pairs of
$\SU_3(q^2)$.

\subsubsection{Non-degenerate unitary space} \label{nondegspace}

Consider $G\cong\SU_{n+1}(q^2)$ as a matrix group with respect to an orthonormal basis of its natural module and let $\A$ be the amalgam consisting of the block-diagonal subgroups $\SU_2(q^2)$ and $\SU_3(q^2)$ and  $\SU_2(q^2)\times
\SU_2(q^2)$. 
One has to prove that the universal enveloping group of the amalgam $\A$ coincides with $G$.  A
natural way to show this is via Tits' Lemma, cf.\ Section \ref{tl}, once one knows a geometry with $G$ as a sufficiently transitive group of automorphisms such that $\cA$ is related to the
amalgam of maximal parabolics induced by the action of $G$. 

Such a geometry $\G_{A_n}$ has been identified in \cite{Aschbacher}, \cite{Aschbacher:1993}, \cite{Das:1994} to be an
$(n+1)$-dimensional non-degenerate unitary space $V$ over $\F_{q^2}$.  The elements of $\G_{A_n}$
are the non-trivial proper non-degenerate subspaces $U$ of $V$, the type of a space $U$ being
its dimension, incidence being defined by symmetrised containment. Using standard terminology from incidence geometry, one-dimensional elements of $\G_{A_n}$ are called {\em points}, two-dimensional elements {\em lines}. Fixing an
orthonormal basis $e_1$, \dots, $e_{n+1}$ of $V$, we consider the action of $G$ as a matrix group on $\G_{A_n}$ with respect to that basis. By Witt's Theorem, see \cite{Scharlau:1995}, this action is flag-transitive, so that we can choose an arbitrary flag $F$ in order to describe the amalgam of parabolics.

This amalgam $\A(\G_{A_n},G,F)$ of parabolics, cf.\ Section \ref{ft}, turns out to have the same universal enveloping group as the amalgam $\A$ consisting of the block-diagonal subgroups $\SU_2(q^2)$ and $\SU_3(q^2)$ and  $\SU_2(q^2)\times
\SU_2(q^2)$ of $G$ by the Redundancy Theorem from Section \ref{redundancy} and by \cite[Lemma 29.3]{Gorenstein/Lyons/Solomon:1995}.

\subsubsection{Decomposing cycles}

The crucial observation for applying Tits' Lemma (Section \ref{tl}) and the Redundancy Theorem (Section \ref{redundancy}) is that $\G_{A_n}$ is almost always simply connected and has many simply connected residues.
In \cite{Bennett/Shpectorov:2004} this simple connectedness is shown by proving that every
cycle of the flag complex of $\G_{A_n}$ is null-homotopic, while in \cite{Das:1994} it is proved in odd characteristic by studying certain subgroup complexes of $\SU_{n+1}(q^2)$. 

In this survey I will sketch the proof given in \cite{Bennett/Shpectorov:2004}. Fixing the base
element $x$ to be a point, a standard
technique based on residual connectedness allows to reduce every cycle of $\G_{A_n}$ to a cycle in the
point-line incidence graph, i.e., the graph on the elements of dimension one and two with incidence as adjacency. 
Furthermore, every cycle in the point-line incidence graph can be understood as a
cycle in the {\em collinearity graph} $\Gamma$ of $\G_{A_n}$, i.e., the graph consisting of the points of $\G_{A_n}$ as vertices in which two vertices are adjacent if and only if they lie on a common line of $\G_{A_n}$. 
A cycle in $\Gamma$ that is contained entirely within the residue of an
element of $\G_{A_n}$ is called {\em geometric} and, being contained in a cone, is
null-homotopic.  Thus, simple connectedness of $\G_{A_n}$ follows, if one can prove that every cycle in $\Gamma$ can be decomposed into a product of
geometric cycles.

A key fact exploited in \cite{Bennett/Shpectorov:2004} is that up to a few exceptions
$\Gamma$ has diameter two.  This implies that every cycle in
$\Gamma$ is a product of cycles of length up to five and, thus, it suffices to
show that every cycle of length three, four, and five is null-homotopic.
When the dimension is large, one can always find a point that is perpendicular
to all points on a fixed cycle, producing a decomposition of that cycle
into geometric triangles. Hence proving simple connectedness is more or less trivial for
large dimension.  The difficulty of the proof lies in the case of small
dimension, where \cite{Bennett/Shpectorov:2004} resorts to a case-by-case analysis. 

To give the precise statement, let $n \geq 3$ and let $q$ be any prime power. Then the geometry $\G_{A_n}$ is simply connected, if $(n,q)$ is not
one of $(3,2)$ and $(3,3)$. Since neither of these
exceptions is simply connected, cf.\ Section \ref{su43}, the result in \cite{Bennett/Shpectorov:2004} is optimal.

\subsubsection{Phan-type theorem of type $A_n$}

Altogether the Phan-type theorem of type $A_n$ follows:

\begin{phantype}[Bennett, Shpectorov \cite{Bennett/Shpectorov:2004}] \label{main1} \label{main2}Let $q$ be a prime power, let $n \geq 3$, and let $G$ be a group admitting a weak Phan system of type $A_n$ over $\mathbb{F}_{q^2}$.
\begin{enumerate}
\item If $q \geq 4$, then $G$ is isomorphic to a central quotient of $\SU_{n+1}(q^2)$.
\item If $q = 2$, $3$ and $n \geq 4$ and if, furthermore,
\begin{enumerate}
\item for any triple $i$, $j$, $k$ of nodes of the Dynkin diagram $A_n$ that form a subdiagram $$\node^{i}\stroke{}\node^{j}\stroke{}\node^{k}$$ of type $A_3$, the subgroup 
$\gen{U_{i,j},U_{j,k}}$ is isomorphic to a central quotient of 
$\SU_4(q^2)$;
\item in case $q=2$
\begin{enumerate}
\item for any triple $i$, $j$, $k$ of nodes of the Dynkin diagram $A_n$ that form a subdiagram $$\node^{i} \quad \quad \node^{j}\stroke{}\node^{k}$$ of type $A_1 \oplus A_2$ the groups $U_i$ and $U_{j,k}$ commute 
elementwise; and 
\item for any quadruple $i$, $j$, $k$, $l$ of nodes of the Dynkin diagram $\Delta$ that form a subdiagram $$\node^{i}\stroke{}\node^{j} \quad \quad \node^{k}\stroke{}\node^{l}$$ of type $A_2 \oplus A_2$
the groups
$U_{i,j}$ and $U_{k,l}$ commute elementwise;
\end{enumerate}
\end{enumerate}
then $G$ is isomorphic to a central quotient of $\SU_{n+1}(q^2)$.
\end{enumerate}
\end{phantype}

\subsubsection{The group $\SU_4(3^2)$} \label{su43}

The extra conditions in Main Theorem \ref{main2} (ii) are due to the fact that for small $q$ and $n$ the geometry $\G_{A_n}$ is not simply connected.
For example, \cite{Horn:2008a} describes a group $H$ admitting a weak Phan system of type $A_3$ over $\mathbb{F}_{3^2}$ that is isomorphic to a non-split central extension of $\SU_4(3^2)$ by a group $K \cong (\mathbb{Z}/3\mathbb{Z})^2$, i.e. the sequence $1\to K \to H \rightarrow \SU_4(3^2) \to 1$ is exact and non-split; in fact, $H$ is isomorphic to the Schur cover of $\SU_4(3^2)$. From there it is deduced in \cite{Horn:2008a} that the geometry $\G_{A_3}$ admits a $9$-fold universal cover in case $q=3$. 

\section{The Curtis--Tits Theorem} \label{Curtis--Tits theorem}

Phan's theorems can be considered as a twisted version of the Curtis--Tits Theorem. Therefore by explaining the general setup of Phan-type theorems one naturally also describes a setup of the Curtis--Tits Theorem. In this section I will give many different (sometimes inequivalent) ways how to state the Curtis--Tits Theorem. Some versions deal with determining a Chevalley group (or even a Kac--Moody group) as the universal enveloping group of a certain amalgam, others with characterisations of these groups from purely local data. One version is merely concerned with the simple connectedness of a suitable chamber system. Each version has its advantages and disadvantages. While it may be easier for the geometric group theorist to prove the simple connectedness of some complex, a local group theorist may prefer to apply a version requiring only knowledge about local data. The transition from the former point of view to the latter requires a certain amount of rigidity of the complex on which the group of interest acts. This can be exploited to obtain a classification of amalgams as achieved in \cite{Bennett/Shpectorov:2004}, \cite{Dunlap:2005} or, more ambitiously, to obtain a classification of groups generated by a class of abstract root groups as sketched in Section \ref{timm}.   

\subsection{The result} \label{ct}

\subsubsection{Chevalley groups and the Steinberg presentation} \label{cgsp}

Chevalley groups can be defined by their Steinberg presentation, cf.\ \cite[Theorem 8]{Steinberg:1968}, the approach I decided to take in this survey. For additional background information and terminology see \cite{Bourbaki:2002}, \cite{Gorenstein/Lyons/Solomon:1998}, \cite{Steinberg:1962}, \cite{Steinberg:1968}.

Let $\Sigma$ be an indecomposable root system of rank at least two and let $\F$ be a field. Consider the group $G$ generated by the collection of elements $$\lb x_r(t) \mid r \in \Sigma, t\in \F \rb$$ subject to the following relations:
\begin{enumerate}
\item $x_r(t)$ is additive in $t$.
\item If $r$ and $s$ are roots and $r + s \neq 0$, then
$$\lbr x_{r}(t),x_s(u) \rbr = \prod x_{hr+ks}(C_{h k r s}t^hu^k)$$
with $h, k > 0$, $hr+ks \in \Sigma$ (if there are no such numbers, then $\lbr x_r(t), x_s(u) \rbr = 1$), and certain structure constants $C_{h k r s} \in \lb \pm 1, \pm 2, \pm 3 \rb$.
\item $h_r(t)$ is multiplicative in $t$, where $h_r(t)$ equals $w_r(t)w_r(-1)$ and $w_r(t)$ equals $x_r(t)x_{-r}(-t^{-1})x_r(t)$ for $t \in \F^*$.
\end{enumerate}
With the correct choice of the structure constants $C_{h k r s}$ (see \cite[Theorem 1.12.1]{Gorenstein/Lyons/Solomon:1998}, \cite{Steinberg:1968}) the group $G$ is called the {\em universal Steinberg--Chevalley group} constructed from $\Sigma$ and $\F$. For $r \in \Sigma$ the group $$x_r = \lb x_r(t) \mid t \in \F \rb = (\F, +),$$ and any conjugate of $x_r$ in $G$, is called a {\em root (sub)group}.
By \cite[Theorem 9]{Steinberg:1968}, if $\Sigma$ is an indecomposable root system of rank at least two and $\F$ an algebraic extension of a finite field, then the above relations {\rm (i)} and {\rm (ii)} suffice to define the corresponding universal Chevalley group, i.e., they imply the relations {\rm (iii)}.

\subsubsection{Redundancy of the Steinberg presentation} \label{redundancy2}

The Curtis--Tits Theorem states that Steinberg's presentation of Chevalley groups in Section \ref{cgsp} is highly redundant and that the amalgam consisting of rank one and rank two subgroups with respect to a system of fundamental roots of a maximal torus of a Chevalley group suffices to present this Chevalley group, cf.\ \cite{Curtis:1965}, \cite{Tits:1962}, \cite[Theorem 13.32]{Tits:1974}. 

The following version of the Curtis--Tits Theorem refers directly to the Steinberg presentation.

\begin{curtistits}[{Curtis \cite[Corollary 1.8]{Curtis:1965}}] \label{Curtis}
Let $\Sigma$ be an indecomposable root system of rank at least two, let
$\Pi$ be a fundamental system of $\Sigma$, and let $\F$ be an arbitrary field with five distinct elements. Define $G$ to be the abstract group with 
generators $\lb x_r(t) \mid r \in \Sigma, t \in \F \rb$ and defining 
relations
\begin{eqnarray}
x_r(t)x_r(u) = x_r(t+u), r \in \Sigma, t, u \in \F,   \label{first}
\end{eqnarray}
and for independent roots $r$, $s$,
\begin{eqnarray}
\lbr x_r(t), x_s(u) \rbr = \prod x_{hr+ks} (C_{h k r s} t^h u^k), \label{second}
\end{eqnarray}
with $h, k > 0$, $hr+ks \in \Sigma$ (if there are no such numbers, then $\lbr x_r(t), x_s(u) \rbr = 1$), and structure constants $C_{h k r s} \in \lb \pm 1, \pm 2, \pm 3 \rb$. 

Let $A = \bigcup A_{ij}$, where 
$A_{ij}$ is the set of all roots which are linear combinations of the 
fundamental roots $r_i, r_j \in \Pi$. Let $G^*$ be the abstract group with 
generators $\lb x_r(t) \mid r \in \Sigma, t \in \F \rb$ and defining 
relations {\rm (\ref{first})}, for $r \in A$, and {\rm (\ref{second})} for independent roots $r$, $s$ belonging to some $A_{ij}$. 

Then the natural epimorphism $G^* \rightarrow G$ is an isomorphism.
\end{curtistits}

A more compact formulation (albeit without a concrete presentation) can be found in \cite{Gorenstein/Lyons/Solomon:1998}, \cite{Tits:1962}, \cite[Theorem 13.32]{Tits:1974}. Generalisations and variations on the theme are contained in \cite{Caprace:2007a}, \cite{Timmesfeld:2004}.

\begin{curtistits}[Gorenstein, Lyons, Solomon \cite{Gorenstein/Lyons/Solomon:1998}, Tits \cite{Tits:1974}] \label{tits}
Let $G$ be the universal version of a Chevalley group of
(twisted) rank at least three with root system $\Sg$, fundamental system
$\Pi$, and root groups $X_\alpha$, $\alpha \in \Sg$. For each $J
\subseteq \Pi$ let $G_J$ be the subgroup of $G$ generated by all root
subgroups $X_\alpha$, $\pm \alpha \in J$. Let $D$ be the set of all
subsets of $\Pi$ with at most two elements. Then $G$ is the universal
enveloping group of the amalgam $(G_J)_{J \in D}$.
\end{curtistits}

To look at a concrete example, consider the case of the
universal Steinberg--Chevalley group of type $A_n$, which is $G=\SL_{n+1}(\F)$. With the usual
choices of the root subgroups in $G$ and of a basis of the natural module of $G$, the subgroups $G_{\alpha,\beta}$ generated by the {\em fundamental rank one subgroups} $G_\alpha := \langle X_\alpha, X_{-\alpha} \rangle$ and $G_\beta := \langle X_\beta, X_{-\beta} \rangle$ are
the block-diagonal subgroups $\SL_3(\F)$ and $\SL_2(\F)\times \SL_2(\F)$.

\subsubsection{The Curtis--Tits Theorem Phan-style} \label{ctps} \label{uoar}

The Curtis--Tits Theorem has been extended to a result including a classification of amalgams by Phan \cite{Phan:1970} (for $\SL_{n+1}(q)$), by Humphreys \cite{Humphreys:1972} (for every finite Chevalley group with a simply laced diagram), and by Dunlap \cite{Dunlap:2005} (for every Chevalley group).
 Phan constructs a $BN$-pair, cf.\ Section \ref{BN}, from the amalgam he is starting with and consequently recognises his target group as a group with a $BN$-pair of type $A_n$. Humphreys \cite{Humphreys:1972} gives another proof of the main result of \cite{Phan:1970} whose central idea is identical to Bennett and Shpectorov's \cite{Bennett/Shpectorov:2004} proof of uniqueness of Phan amalgams. After obtaining uniqueness Humphreys \cite{Humphreys:1972} simply invokes the Curtis--Tits Theorem. He mentions in passing that Curtis--Tits amalgams can be classified for the non-simply laced spherical diagrams of rank at least three, if one can control the behaviour of the root subgroups of $\Sp_4(q)$. Similarly, Shpectorov mentioned to me that a classification of Phan amalgams can be accomplished as soon as one can control the behaviour of the Phan amalgam in $\Sp_6(q)$. Both observations have been worked out in detail by now, see \cite{Dunlap:2005}, \cite{Gramlich:2004d}. 

I point out here that Timmesfeld has also obtained proofs of the Curtis--Tits Theorem. One approach is also based on the construction of $BN$-pairs, see \cite{Timmesfeld:2004}, while an alternative approach (see \cite{Timmesfeld:1998}, \cite{Timmesfeld:2003}, \cite{Timmesfeld:2006}) is based on his theory of abstract root subgroups \cite{Timmesfeld:1991}, \cite{Timmesfeld:1999}, \cite{Timmesfeld:2001}; see Section \ref{timm}.

\begin{curtistits}[Phan \cite{Phan:1970}, Humphreys \cite{Humphreys:1972}, Timmesfeld \cite{Timmesfeld:2004}, Dunlap \cite{Dunlap:2005}] \label{phanstyle}
Let $\Delta$ be a spherical Dynkin diagram of rank at least three, let $\mathbb{F}$ be a field, and let $G$ be a group generated by subgroups $G_\alpha$ and $G_{ \alpha, \beta }$, for all $\alpha, \beta \in \Delta$, isomorphic to Chevalley groups over $\mathbb{F}$ as indicated by the induced Dynkin diagram on the nodes $\alpha$, $\beta$, with the property that in each $G_{\alpha, \beta}$ the subgroups $G_\alpha$ and $G_\beta$ correspond to the choice of a fundamental system of roots with respect to a maximal torus of $G_{\alpha, \beta}$. Then $G$ is a central quotient of the universal Chevalley group of type $\Delta$ over $\mathbb{F}$. 
\end{curtistits}

Note that a theorem like the Curtis--Tits Theorem Version \ref{phanstyle} is much easier to apply in local group theory than the Curtis--Tits Theorem Versions \ref{Curtis} and \ref{tits}.

\subsection{Buildings and twin buildings} \label{twin buildings}

\subsubsection{Towards a Curtis--Tits geometry}

One purpose of this survey is to point out similarities between the Curtis--Tits Theorem on one hand and Phan's theorems on the other hand by describing suitable geometries whose simple connectedness yields the respective group-theoretic identification result via Tits' Lemma \cite{Tits:1986} (Section \ref{tl}). These geometries can be constructed using the opposition relation of a building or a twin building, cf.\ \cite{Abramenko/Muehlherr:1997}, \cite{Muehlherr}, \cite{Tits:1989}. 
I start with a description of the ideas of proof of the Curtis--Tits Theorem given in \cite{Abramenko/Muehlherr:1997}, \cite{Muehlherr}, whose generality actually implies that result for any two-spherical diagram (i.e., each sub-diagram of cardinality at most two is spherical), except that one has to exclude some small cases covered by the original Curtis--Tits Theorem. These exceptions arise from exactly those rank two diagrams and fields for which the geometry opposite a chamber inside the corresponding Moufang polygon is not connected, see \cite{Abramenko/Maldeghem:1999}, \cite{Brouwer:1993}.
Before I am able to properly explain this geometric approach to the Curtis--Tits Theorem, I need to introduce the concepts of a building, a twin building, and of the opposite geometry.

A Chevalley group $G$ acts on its natural geometry,
called a {building}. 
Buildings have been developed by Tits in numerous articles since the mid-1950's. The standard reference are Tits' lecture notes \cite{Tits:1974}. Other references are \cite{Abramenko/Brown:2008}, \cite{Brown:1989}, \cite{Grundhoefer:2002}, \cite{Ronan:1989}, \cite{Scharlau:1995}, \cite{Weiss:2003}. For Coxeter groups and root systems see \cite{Bourbaki:2002}, \cite{Davis:2008}, \cite{Humphreys:1990}.

\subsubsection{Coxeter systems} \label{cox}

For a {\em Coxeter matrix} $M = (m_{ij})_{i, j \in I}$ over some finite set
$I$, i.e., a symmetric $|I| \times |I|$-matrix over $\mathbb{N} \cup \{ \infty \}$ whose diagonal entries equal to one and whose off-diagonal entries are greater or equal two, the {\em Coxeter diagram} of $M$ is the complete labelled graph with vertex set $I$ and labels $m_{ij}$ on the edge $\{ i, j \}$. The cardinality $|I|$ is called the {\em rank} of the Coxeter diagram. Usually the edges with label $2$ are erased, so that it is meaningful to talk about connected or disconnected Coxeter diagrams.

Let $(W,S)$ be the {\em Coxeter system of type $M$}, i.e., $S = \lb s_i \mid i \in I \rb$ is a set and $W = \gen{ S \mid (s_is_j)^{m_{ij}} = 1}$ is the quotient of the free group generated by $S$ and subject to the relations given by the Coxeter matrix $M$. The Coxeter system is {\em spherical}, if $|W| < \infty$, and {\em irreducible}, if the Coxeter diagram is connected. Irreducible spherical Coxeter diagrams have been classified, cf.\ \cite{Coxeter:1935}. 

Using the Bourbaki notation, the irreducible spherical Coxeter systems of rank at least three fall into the families $A_n$, $B_n$, $C_n$, $D_n$ plus the exceptional diagrams $E_6$, $E_7$, $E_8$, $F_4$. If $\Delta$ is a Coxeter diagram of type $M$, a Coxeter system $(W,S)$ of type $M$ is also called a {\em Coxeter system of type $\Delta$}. For $J \subseteq I$, the pair $(W_J,S_J)$ consisting of $S_J = \{ s_i \in S \mid i \in J \}$ and $W_J = \gen{S_J}$ is also a Coxeter system satisfying $W_J  = \gen{S_J \mid (s_is_j)^{m_{ij}} = 1}$ by \cite[Section IV.1.8, Theorem 2]{Bourbaki:2002}. The group $W$ of a Coxter system $(W,S)$ is called a {\em Coxeter group}. It is in general not possible to reconstruct the Coxeter system from the abstract group $W$, see \cite{Bahls:2005}, \cite{Caprace/Muehlherr:2007}, \cite{Charney/Davis:2000}.  

\subsubsection{Buildings} \label{buildings}

A {\it building of type $(W,S)$} (where $(W,S)$ is a Coxeter system) is a pair $\cB=(\cC, \dl)$
where $\cC$ is a set and the
{\em distance function} $\dl: \cC \times \cC \rightarrow W$  satisfies the following axioms for $x, y \in \cC$ and $w=\dl(x,y)$.
  \begin{description}
        \item[{\bf (B1)}] $w=1$ if and only if $x=y$,

        \item[{\bf (B2)}] if $z \in \cC$ is such that $\dl(y,z)= s \in S$,
        then $\dl(x,z)= w$ or $ws$; furthermore if $l(ws)=l(w)+1$, then
        $\dl(x,z)= ws$, and

        \item[{\bf (B3)}] if $s \in S$, there exists $z \in \cC$ such that
        $\dl(y,z)=s$ and $\dl(x,z)=ws$.
  \end{description}
The group $W$ is called the {\em Weyl group} of the building $\cB$. The building $\cB$ is called {\em spherical}, if its Weyl group $W$ is finite.
Given a building $\cB=(\cC, \dl)$ one can define a chamber system on
$\cC$ in which two chambers $c$ and $d$ are $i$-adjacent, in symbols $c \sim_i d$, if and only if
$\dl(c,d)=s_i$ or $\dl(c,d)=1$. The chamber system  $( \mathcal{C},(\sim_i)_{i \in I})$ 
uniquely determines ${\cal B}$, i.e., the $i$-adjacency relations
on $\mathcal{C}$ determine the distance function $\delta$; cf.\ \cite{Tits:1981}.

In this survey  we only consider buildings
$\cB$
for which the chamber system $\cC$ is thick. All thick spherical buildings with a connected Coxeter diagram $\Delta$ of rank at least three ($|\Delta|$ is also called the {\em rank} of the building) are known, e.g., by a local to global approach using the classification of Moufang buildings of rank two, see \cite[Chapter 40]{Tits/Weiss:2002}, also \cite{Weiss:2003}. This local to global approach is possible, because all thick spherical buildings of rank at least three with a connected Coxeter diagram are Moufang (see \cite[Addendum]{Tits:1974}), whence their rank two residues are Moufang. 
Buildings of rank two are called {\em generalised polygons}, and are studied | Moufang or not | in \cite{Hughes/Piper:1973}, \cite{Payne/Thas:1984}, \cite{Pickert:1975}, \cite{Maldeghem:1998}. 

If $\cB$ is a building, its chamber system
contains a class of thin sub-chamber systems called {\it apartments}, each of which forms a building of the same type as $\cB$.  In an
apartment $\Sg$, for any $c \in \Sg$ and $w \in W$, there is a unique
chamber $d \in \Sg$ such that $\dl(c,d)=w$.  Every pair of chambers
of
$\cC$ is contained in an apartment, cf.\ \cite[Corollary 8.6]{Weiss:2003}. The chamber system
$\cC$ defined by a building is always geometric; indeed buildings have the property, cf.\ Section \ref{csap}, that for any set $\{ (x_i,i) \mid i \in I \}$, with $x_i$ an $I \backslash \{ i \}$-residue of $\mathcal{C}$, such that $x_i \cap x_j \neq \emptyset$ for all $i, j \in I$, the intersection $\bigcap_{i \in I} x_i$ is non-empty, and that for distinct chambers $c$, $d$ of $\mathcal{C}$ there is some $I \backslash \{ i \}$-residue of $\mathcal{C}$ containing $c$ but not $d$, see \cite[Section 1.6]{Tits:1981}. The geometry $\G(\cB) = \G(\mathcal{C})$ is called the {\em building geometry}.

In the language of algebraic groups the following examples for buildings can be given, see \cite[Theorem 5.2]{Tits:1974}, also
\cite[Section 6.8]{Benson:1991} and \cite[Section 2]{Garland:1973}. Starting with a reductive algebraic group $G$ defined over a field $\F$, the {\em Tits building $\G(G,\F)$ of $G$ over $\F$} consists of the simplicial complex whose simplices are indexed by the parabolic $\F$-subgroups of $G$ ordered by the reversed
inclusion relation on the parabolic subgroups. The Steinberg functor and Chevalley--Demazure group schemes, see \cite{Chevalley:1966}, \cite{Demazure:1965}, \cite{Tits:1987}, allow to construct a vast amount of groups yielding a rich supply of buildings.

A key property of buildings is the {\em gate property}, see \cite{Muehlherr:1994}, \cite{Ronan:1989}, and \cite{Weiss:2003}:
For a chamber $c \in C$ and a $J$-residue $$R_J(d) := \{ z \in C \mid \delta(d,z) \in W_J \} \subset \mathcal{C} \quad \quad \quad  \mbox{(cf.\ Section \ref{2.3.1})}$$ there exists a unique chamber $x \in R_J(d)$ such that for all $y \in R_J(d)$ one has
$\delta(c,y) = \delta(c,x) \delta(x,y)$ and, in particular,
$l(c,y) = l(c,x) + l(x,y)$, where $l$ denotes the length function of $W$ with respect to the generating system $S$.
This chamber $x$ is called the {\it projection of $c$ onto $R_J(d)$} and is denoted by $\proj_{R_J(d)} c$.

Any building $\cB$ (and hence its geometry $\G(\cB)$) of rank at least three is simply connected. In fact, more is known about the homotopy type of a building. 

\medskip \noindent {\bf Solomon-Tits Theorem (Solomon, Tits \cite{Solomon:1969}).} {\em 
A spherical Tits building of rank $n$ is homotopy equivalent to a wedge of spheres of dimension $n-1$. A spherical Tits building of rank $n$ over a field of $q$ elements is homotopy equivalent to a wedge of $q^m$ spheres of dimension $n-1$, where $m$ is the number of positive roots.
}

\medskip
The Solomon-Tits Theorem has numerous applications in representation theory. I refer the interested reader to \cite{Humphreys:1987} for an excellent survey and guide to the literature.

\subsubsection{Tits systems} \label{BN}

Let $G$ be a group and $B$, $N$ be subgroups of $G$. The tuple $(G, B, N, S)$ is called a {\em Tits system}, if the following conditions are satisfied:
\begin{enumerate}
\item $G$ is generated by $B$ and $N$;
\item $H = B \cap N$ is normal in $N$;
\item $W = N/H$ admits a finite system $S = \{ w_i \mid i \in I\}$ of generators making $(W,S)$ a Coxeter system;
\item For any $w_i \in S$ we have $w_iBw_i^{-1} \neq B$;
\item For any $w_i \in S$ and all $w \in W$ we have $w_iBw \subseteq (BwB) \cup (Bw_iwB)$.
\end{enumerate}
The pair of subgroups $B$, $N$ of $G$ is also called a {\em $BN$-pair} of $G$, see \cite{Bourbaki:2002}, \cite{Tits:1974}. A group $G$ admitting a $BN$-pair satisfies $$G = \bigsqcup_{i \in I} Bw_iB.$$ For each $i \in I$ the set $P_i := B \cup Bw_iB$ is a subgroup of $G$. A Tits system $(G, B, N, S)$ leads to a building whose set of chambers equals $G/B$ and whose distance function $$\dl : G/B \times G/B \to W$$ is given by $\dl(gB,hB) = w$ if and only if $Bh^{-1}gB = BwB$. In the corresponding chamber system $gB$ and $hB$ are $i$-adjacent if and only if $Bh^{-1}gB \subseteq B \cup Bw_iB$. 

\subsubsection{Twin buildings} \label{twinbuild}

The simple connectedness of a building does not imply the Curtis--Tits Theorem, since the action of a Chevalley group on its building does not yield the correct amalgam.
A class of geometries that yields the correct amalgams is best described using twin buildings.
Twin buildings are obtained by relating two Tits buildings via a co-distance function, see \cite{Abramenko:1996}, \cite{Abramenko/Brown:2008}, \cite{Muehlherr:2002}, \cite{Ronan:2002}, \cite{Ronan/Tits:1994}, \cite{Tits:1992}.
Given two buildings $\cB_+=(\cC_+, \dl_+)$ and $\cB_-=(\cC_-, \dl_-)$ of
the same type $(W,S)$, a {\it co-distance}, also called {\it twinning}, is a map
$$\dl_*:
(\cC_+\times \cC_-)\cup (\cC_-\times \cC_+) \rightarrow W$$ such
that the following axioms hold where $\eps = \pm$, $x \in
\cC_{\eps}$, $y \in \cC_{-\eps}$, $w=\dl_*(x,y)$:
\begin{description}
        \item[{\bf (T1)}] $\dl_*(y,x)=w^{-1}$,

        \item[{\bf (T2)}] if $z \in \cC_{-\eps}$ with $\dl_{-\eps}(y,z)=s\in
S$
        and $l(ws)=l(w)-1$, then $\dl_*(x,z)=ws$, and

        \item[{\bf (T3)}] if $s\in S$, there exists $z \in \cC_{-\eps}$ such that
        $\dl_{-\eps}(y,z)=s\in S$ and $\dl_*(x,z)=ws$.
\end{description}
A {\it twin building} of type $(W,S)$ is a triple $(\cB_+,\cB_-,\dl_*)$,
where $\cB_+$ and $\cB_-$ are buildings of type $(W,S)$ and $\dl_*$ is a
twinning between $\cB_+$ and $\cB_-$. 

Every spherical twin building can be obtained in a unique way from some building $\cB=(\cC,\dl)$ of the same type $(W,S)$, cf.\  \cite[Proposition 1]{Tits:1992}. Let
$\cB_+=(\cC_+,\dl_+)$ be a copy of $\cB$, let $\cB_-=(\cC_-,\dl_-)$
be $(\cC,w_0 \dl w_0)$, and let $\dl_*$ be $w_0\dl$ on $\cC_+ \times \cC_-$ and $\dl
w_0$ on $\cC_- \times \cC_+$, where $w_0$ is the longest element of the Weyl group $W$.

If $R$ is an arbitrary spherical residue of type $J$ in a twin building, then by \cite[4.1]{Ronan:2000} there is a unique 
chamber $z \in R$ with 
$(\delta_*(c,z))_{W_J} = \delta_*(c,z)$ in analogy to the gate property of a building. Moreover, by \cite[4.3]{Ronan:2000}, for all $y \in R$ we have
$\delta_*(c,y) = \delta_*(c,z) \delta_-(z,y)$ and in particular
$l_*(c,y) = l_*(c,z)-l(z,y)$. As for buildings,
this chamber $z$ is called the {\it projection} of $c$
onto $R$ and is denoted by $\proj_R c$.
Furthermore, if $J$ is a spherical subset of $S$, then any two $J$-residues
of ${\cal B}_{\epsilon}$ are isomorphic for each $\epsilon \in \{ +,- \}$.
Additionally, there exists a twin version of the
main result in \cite{Dress/Scharlau:1987}, as observed in \cite{Devillers/Muehlherr:2007}, stating that, if $R$, $Q$ are spherical residues of a twin building, then $\proj_R Q := \{ \proj_R x \mid x \in Q \}$ is a spherical residue contained in $R$. Moreover, for $ R' := \proj_R Q$ and $Q' := \proj_Q R$, the maps
$\proj_{R'}^{Q'} := {\proj_{R'}}_{|Q'}: Q' \rightarrow R'$ and 
$\proj_{Q'}^{R'} := {\proj_{Q'}}_{|R'}: R' \rightarrow Q'$ are adjacency-preserving
bijections inverse to each other.

\subsection{The opposite geometry and its simple connectedness} \label{opp}

\subsubsection{Opposition}

The concept of the opposite geometry can be traced to Tits \cite{Tits:1989}. The opposition relation is an important concept in the theory of buildings and plays a crucial role in \cite{Abramenko:1996}, \cite{Abramenko/Ronan:1998}, \cite{Abramenko/Maldeghem:1999}, \cite{Abramenko/Maldeghem:2000}, \cite{Abramenko/Maldeghem:2001}, \cite{Abramenko/Maldeghem:2002}, \cite{Muehlherr}, \cite{Muehlherr:1998}. Given a twin building $\cT=(\cB_+,\cB_-,\dl_*)$, one can define the
chamber system $\Opp(\cT)$ on the set $$\{(c_+,c_-)\in \cC_+\times \cC_- \mid\,
\dl_*(c_+,c_-)=1 \}$$ in which $(c_+,c_-) \sim_i (d_+,d_-)$ if and only if $c_+ \sim_i d_+$ and $c_- \sim_i d_-$.  Chambers $x\in\cC_+$ and $y\in\cC_-$ with
$\dl_*(x,y)=1$ are called {\it opposite}, hence the notation. 

Denote the corresponding pregeometry by $\G_{{\rm op}}$. For $\G_+$ and
$\G_-$ the building geometries that correspond to $\cB_+$ and
$\cB_-$,  elements $x_+\in\G_+$ and $x_-\in\G_-$ of the same type
$i\in I$ are called {\em opposite}, if they are contained in opposite
 chambers.  The elements of the pregeometry $\G_{{\rm op}}$ of type $i$
are the pairs $(x_+,x_-)$ of opposite elements of type $i$. Two pairs
$(x_+,x_-)$ and $(x'_+,x'_-)$ are incident in $\G_{{\rm op}}$, if $x_+$
and $x'_+$ are incident in $\G_+$ and $x_-$ and $x'_-$ are incident
in $\G_-$.  Clearly, a pair $(c_+,c_-)\in\Opp(\cT)$ produces a
maximal flag in $\G_{{\rm op}}$, and it can be shown that every maximal
flag is obtained in this way.
Hence the pregeometry $\G_{{\rm op}}$ is a geometry, called the {\it
opposite geometry}. Moreover, the chamber system $\Opp(\cT)$ is geometric which follows by a building-theoretic argument proving that the map $c \mapsto \psi_\cC(c)$ in Section \ref{chsys} is an isomorphism. 

\subsubsection{Examples of classical opposite geometries} \label{exclass}

The following examples are descriptions of the opposite geometries for the four classical series of spherical buildings.

\medskip \noindent
{\bf Example 1a.} \label{example1a} Let $\F$ be an arbitrary field and consider the universal Steinberg--Chevalley group $G\cong \SL_{n+1}(\F)$ of type $A_n$ over $\F$.
It corresponds to the building geometry $\G$ of type $A_n$, better known as the projective
geometry, whose elements of type $i$, $1\le i\le n$, are the
$i$-dimensional subspaces in an $(n+1)$-dimensional $\F$-vector
space $V$.  The geometries $\G_+$ and $\G_-$ are isomorphic to the projective geometry $\G$ and its dual, respectively. The latter is identical to $\G$
except that the types are interchanged by the map $i \mapsto n+1-i$.
Elements $x_+\in \G_+$ and $x_-\in \G_-$ of type $i$ are
opposite if they intersect trivially or, equivalently, form a direct sum
decomposition $V=x_+\oplus x_-$, cf.\ \cite[II, \S4, Lemma 23]{Abramenko:1996}.  These decompositions
are the elements of $\G_{{\rm op}}$, where $x_+\oplus x_-$ is incident to $x'_+\oplus x'_-$ if and only if $x_\epsilon$ is incident to $x'_\epsilon$ for $\epsilon \in \{ +, - \}$.

\medskip \noindent
{\bf Example 2a.} Let $G\cong \Spin_{2n+1}(\F)$ be the universal Steinberg--Chevalley group corresponding to the building geometry $\G$ of
type $B_n$.  The geometry $\G$ is the geometry of all totally isotropic
subspaces of a non-degenerate $(2n+1)$-dimensional orthogonal space $V$ over $\F$.
In this case, both $\G_+$ and $\G_-$ are isomorphic to $\G$ and two
$i$-dimensional
totally isotropic subspaces $x_+$ and $x_-$ are opposite if $x_-$
intersects the orthogonal complement of $x_+$ trivially, i.e., $x_+^\perp \cap x_- = \{ 0 \}$ or, equivalently, $x_+^\perp \oplus x_- = V$, \cite[II, \S6, Lemma 29]{Abramenko:1996}.  Such
pairs $(x_+, x_-)$ are the elements of $\G_{{\rm op}}$, where $(x_+,x_-)$ is incident to $(x'_+,x'_-)$ if and only if $x_\epsilon$ is incident to $x'_\epsilon$ for $\epsilon \in \{ +, - \}$.

\medskip \noindent
{\bf Example 3a.} Consider the universal Steinberg--Chevalley group $G\cong \Sp_{2n}(\F)$ of type $C_n$. In this case the corresponding building geometry $\G$ is the geometry of all totally isotropic
subspaces of a non-degenerate $2n$-dimensional symplectic space $V$ over $\F$.
Both $\G_+$ and $\G_-$ are isomorphic to $\G$.  Two
$i$-dimensional
totally isotropic subspaces $x_+$ and $x_-$ again are opposite if $x_-$
intersects the orthogonal complement of $x_+$ trivially, i.e., $x_+^\perp \cap x_- = \{ 0 \}$ or, equivalently, $x_+^\perp \oplus x_- = V$, \cite[II, \S6, Lemma 29]{Abramenko:1996}. The
pairs $(x_+, x_-)$ are the elements of $\G_{{\rm op}}$, where $(x_+,x_-)$ is incident to $(x'_+,x'_-)$ if and only if $x_\epsilon$ is incident to $x'_\epsilon$ for $\epsilon \in \{ +, - \}$.

\medskip \noindent
{\bf Example 4a.} Let $G\cong \Spin^+_{2n}(\F)$ be the universal Steinberg--Chevalley group of type $D_n$, to which corresponds the building geometry $\G$ of totally isotropic
subspaces of a non-degenerate $2n$-dimensional orthogonal space $V$ over $\F$ of Witt index $n$. In this case, both $\G_+$ and $\G_-$ are isomorphic to $\G$ up to interchanging the elements of types $n-1$ and $n$ in case $n$ odd.  Two
totally isotropic subspaces $x_+$ and $x_-$ of type $i$ are opposite, if $x_-$
intersects the orthogonal complement of $x_+$ trivially, i.e., $x_+^\perp \cap x_- = \{ 0 \}$ or, equivalently, $x_+^\perp \oplus x_- = V$, \cite[II, \S7, Lemma 31]{Abramenko:1996}.  Such
pairs $(x_+, x_-)$ are the elements of $\G_{{\rm op}}$, where $(x_+,x_-)$ is incident to $(x'_+,x'_-)$ if and only if $x_\epsilon$ is incident to $x'_\epsilon$ for $\epsilon \in \{ +, - \}$.

\subsubsection{The Curtis--Tits Theorem via geometric group theory}

If a twin building admits a strongly transitive group of automorphisms, i.e., a group acting transitively on the pairs of opposite chambers, the group acts flag-transitively on $\G_{{\rm op}}$. In case the acting group is semisimple split algebraic or Kac--Moody, the stabilisers of the elements of a
maximal flag of $\G_{{\rm op}}$ are products of the type $G_{\Pi\setminus
\{\alpha\}} Z(T)$, where $T$ is a maximal torus and $\Pi$ is a system of fundamental roots with respect to $T$. This setup together with Tits' Lemma (Section \ref{tl}) and the Redundancy Theorem \cite[Theorem 3.3]{Gramlich/Maldeghem:2006} (Section \ref{redundancy}) implies that the Curtis--Tits Theorem stated as in \cite{Gorenstein/Lyons/Solomon:1998} (Section \ref{redundancy2}) follows from the following simple-connectedness result.

\begin{curtistits}[Abramenko, M\"uhlherr \cite{Abramenko/Muehlherr:1997}, M\"uhlherr \cite{Muehlherr}] \label{mueh}
Let $\cT$ be a thick twin building with two-spherical diagram of rank at least three such that there is no rank two residue in $\cB_+$ or $\cB_-$ which is isomorphic to the buildings associated to $B_2(2)$, $G_2(2)$, $G_2(3)$ or $^2F_4(2)$. Then $\Opp(\cT)$ is simply connected.
\end{curtistits}

The proof of this theorem in \cite{Muehlherr} is derived directly
from the axioms of twin buildings, properties of apartments in
buildings, and certain connectedness properties of buildings like their simple connectedness.  The exceptions in this approach come from the fact that the geometry opposite to a chamber in an arbitrary Moufang polygon is connected except in the cases $B_2(2)$, $G_2(2)$, $G_2(3)$ or $^2F_4(2)$, cf.\ \cite{Abramenko/Maldeghem:1999}, \cite{Brouwer:1993}. Note in passing that in \cite[II, \S2, Proposition 9]{Abramenko:1996} it is shown that there is no hope for general connectedness results in the non-Moufang case. Of course, the Curtis--Tits Theorem for Steinberg--Chevalley groups does not have any exceptions by \cite{Curtis:1965}, \cite{Gorenstein/Lyons/Solomon:1998}, \cite{Timmesfeld:2003}, \cite{Tits:1962}, \cite{Tits:1974}.

In \cite{Abramenko/Muehlherr:1997} the logic of proof is turned around. The authors prove the combinatorial Curtis--Tits Theorem Version \ref{mueh} by directly proving the following generalisation of \cite[Theorem 13.32]{Tits:1974}. The key is to construct an RGD system for $G$, cf.\ \cite{Tits:1992}, also \cite[1.5]{Remy:2002}, \cite[Chapter 8]{Abramenko/Brown:2008}.

\begin{curtistits}[Abramenko, M\"uhlherr \cite{Abramenko/Muehlherr:1997}, M\"uhlherr \cite{Muehlherr}]
Let $\cT$ be a thick twin building with two-spherical diagram $\Delta$ of rank at least three such that there is no rank two residue in $\cB_+$ or $\cB_-$ which is isomorphic to the buildings associated to $B_2(2)$, $G_2(2)$, $G_2(3)$ or $^2F_4(2)$, let $G$ be a group acting transitively on the pairs of opposite chambers of $\cT$, and let $(c_+,c_-)$ be a pair of opposite chambers in $\cT$. For each $J
\subseteq \Delta$ let $G_J$ be the subgroup of $G$ stabilizing the $J$-residue of $c_+$ and the $J$-residue of $c_-$. Let $D$ be the set of all
subsets of $\Delta$ with at most two elements. Then $G$ is the universal
enveloping group of the amalgam $(G_J)_{J \in D}$.
\end{curtistits}

A variation on this theme can also be found in \cite{Caprace:2007a}. The classification of amalgams in \cite{Dunlap:2005} allows to formulate this Curtis--Tits Theorem in Phan-style, cf.\ \cite{Humphreys:1972}, \cite{Phan:1970}, also Section \ref{ctps}. 

\subsection{Abstract root subgroups} \label{timm}

A completely different and independent 
approach to the Curtis--Tits Theorem is based on the classification of groups generated by a class of abstract root subgroups \cite{Timmesfeld:1991}, \cite{Timmesfeld:1999}, \cite{Timmesfeld:2001}.
This wonderful classification result makes it possible to prove all sorts of generalisations of Steinberg-presentation-type results and the Curtis--Tits Theorem, see \cite{Timmesfeld:1998}, \cite{Timmesfeld:2003}, \cite{Timmesfeld:2006}. The case of simply laced diagrams is stated more easily than the general case, cf.\  \cite{Timmesfeld:1998}. Hence I restrict myself to presenting that case here. The article \cite{Timmesfeld:2003} deals with every spherical diagram.

\begin{curtistits}[Timmesfeld \cite{Timmesfeld:1998}] \label{timmesfeld}
Let $\Delta$ be a spherical simply laced diagram of rank at least three and let $G$ be a group generated by subgroups $X_i$, $i \in \Delta$, satisfying
\begin{enumerate}
\item $X_i$ is a perfect central extension of $\PSL_2(\F)$, $\F$ a division ring,
\item in each $X_i$ there exists a non-trivial diagonal subgroup $H_i$ normalising all $X_j$, $j \in \Delta$,
\item for $i \neq j$ one of the following holds:
\begin{enumerate}
\item $[X_i,X_j] = 1$;
\item for $X_{ij} := \gen{X_i,X_j}$, the quotient $X_{ij} / Z(X_{ij})$ is isomorphic to $\PSL_3(\F)$, where $Z(X_{ij}) \subseteq X_{ij}'$; moreover, the unipotent subgroups of $X_i$, $X_j$ are mapped onto elation subgroups, corresponding to point-line pairs, of $X_{ij} / Z(X_{ij})$. 
\end{enumerate}
\end{enumerate}
Suppose further that, if $|\F| = 4$, then $|Z(X_{ij})| < 12$ for some connected pair $i$, $j$ of nodes of $\Delta$.

Then $G$ is a perfect central extension of $\PSL_{n+1}(\F)$ ($\F$ a division ring), $\PSOmega_{2n}(\F)$, or the adjoint Steinberg--Chevalley group $E_n(\F)$ ($\F$ a commutative field) and there exists a homomorphism mapping the $X_i$ onto the fundamental subgroups.
Furthermore, if each $X_i$ is a factor group of $\SL_2(\F)$, $\F$ a commutative field, then $G$ is a factor group of the universal Steinberg--Chevalley group of type $A_n$, $D_n$, or $E_n$ over $\F$.
\end{curtistits}

I would like to point out that the paper \cite{Timmesfeld:2004} contains another proof of the Curtis--Tits Theorem, one that is independent of the classification of groups generated by a class of abstract root subgroups. Instead it relies on a construction of $BN$-pairs and can be considered as a direct generalisation of \cite{Phan:1970}. For a generalisation to Kac--Moody groups see \cite{Caprace:2007a}. 

\section{Phan-type theorems for finite Chevalley groups} \label{flipflop geometries}

\subsection{From Aschbacher's geometry to the general construction} \label{flipgen}

In this section we discuss how Aschbacher's geometry \cite{Aschbacher} and its simple connectedness initiated Phan theory.

\subsubsection{Non-degenerate unitary space, revisited} \label{nondegspacerev}

{\bf Example 1b.} Consider the situation of Example 1a from Section \ref{exclass}, but change
the field of definition to $\F_{q^2}$, so that $G \cong \SL_{n+1}(q^2)$.
Consider a unitary polarity $\tau$, that is, an involutory isomorphism
from $\G$ onto the dual of $\G$ which is
defined by a non-degenerate hermitian form $\Phi$ on $V$. The map
$\tau$ sends every subspace of $V$ to its orthogonal complement with
respect to $\Phi$ and produces an involutory automorphism of
the twin building $\cT$ that switches $\cC_+$ and $\cC_-$ and, thus,
$\G_+$ and $\G_-$.  It is an automorphism in the sense that it
transforms $\dl_+$ into $\dl_-$ (and vice versa), and preserves
$\dl_*$.  Corresponding to $\tau$, there is an automorphism of $G$, which
is also denoted by $\tau$. The group $$G_\tau=C_G(\tau)\cong
\SU_{n+1}(q^2)$$ acts on
$$\G_\tau=\{(x_+,x_-)\in\G_{{\rm op}} \mid x_+^\tau=x_-\}.$$   The elements of
$\G_\tau$ are of the form $(x_+,x_-)$ where $x_-=x_+^\tau=x_+^\perp$
and $V=x_+\oplus x_-=x_+\oplus x_+^\perp$, cf.\ Example 1a in Section \ref{exclass}.  Thus, the mapping
$(x_+,x_-)\mapsto x_+$ establishes an isomorphism between $\G_\tau$
and the geometry $\G_{A_n}$ of all proper non-degenerate subspaces of the unitary
space $(V, \Phi)$.  This geometry $\G_{A_n}$ is exactly Aschbacher's geometry from Section \ref{nondegspace}.

\subsubsection{Flips and Phan involutions} \label{flipphan}

Section \ref{nondegspacerev} suggests the following general construction introduced in \cite{Bennett/Gramlich/Hoffman/Shpectorov:2003}.  Let $\cT=(\cB_+,
\cB_-,\dl_*)$ be a twin building as defined in Section \ref{twinbuild}. Then an involutory 
automorphism $\tau$ of
$\cT$ satisfying
\begin{description}
\item[{\bf (F1)}] $\cC_+^\tau = \cC_-$,
\item[{\bf (F2)}] $\tau$ flips the distances, \ie
$\dl_\epsilon(x,y)= \dl_{-\epsilon}(x^\tau, y^\tau)$ for $\epsilon = \pm$, and
\item[{\bf (F3)}] $\tau$ preserves the co-distance, \ie $\dl_* (x,y) =
\dl_\ast (x^\tau, y^\tau)$
\end{description}
is called a
{\em flip}. Notice that by (T1) of Section \ref{twinbuild} the element $\delta_*(x,x^\tau)$ always is an involution.

\medskip \noindent
A flip satisfying the additional condition
\begin{description}
\item[{\bf (F4)}] there exists a chamber $c
\in \cC_{\pm}$
with $\dl_*(c,c^\tau)=1$
\end{description}
is called a {\em Phan involution}. 

\medskip \noindent
In case $\tau$ is a Phan involution
the chamber system $\cC_\tau$ whose chambers are pairs
$(c,c^\tau)$ that belong to $\Opp(\cT)$, i.e., $$\cC_\tau = \lb (c_+,c_-) \in \Opp(\mc{T}) \mid \lb c_+, c_- \rb = \lb c_+^\tau, c_-^\tau \rb \rb,$$ is called the {\em flipflop system} of $\tau$.  By (F4) the chamber system $\cC_\tau$ is non-empty. By \cite{Gramlich/Horn/Muehlherr}, \cite{Horn:2008b} very many $\cC_\tau$ are geometric, in particular all flipflop systems encountered in this survey. However, it is not known to me, if $\cC_\tau$ is
geometric in general. For a geometric flipflop system $\cC_\tau$ denote by $\G_\tau$ the
corresponding geometry, the {\em flipflop
geometry}. 

Following \cite{Devillers/Muehlherr:2007} one can alternatively define a Phan involution to be a flip of a twin building satisfying
\begin{description}
\item[{\bf (F4)'}] $\proj_P \tau \neq P$
for each panel $P$ of ${\cal T}$
\end{description}
where $\proj_R \tau :=  \{ x \in R \mid \proj_R \tau(x) = x \}$ for a spherical residue $R$ of ${\cal T}$ (Section \ref{twinbuild}).
It is easily seen that a flip satisfying (F4)' also satisfies (F4). When talking about Phan involutions, we will generally only assume the validity of axioms (F1), (F2), (F3), (F4), unless explicitly stated otherwise.

\subsubsection{Flips of spherical twin buildings} \label{fg}

 For a spherical twin building one can compute the action of
$\tau$ on the Dynkin diagram of the building, see \cite[Section 3.3]{Gramlich:2002}.  Indeed, by Tits' characterisation each spherical twin
building arises from a spherical building $\cB=(\cC,\dl)$ (cf.\ \cite[Proposition 1]{Tits:1992} and also Section \ref{twinbuild} of this survey) and we have $$\dl(c,d)=
\dl_{+}(c,d)=\dl_{-}(c^{\tau},d^{\tau})=w_{0}\dl(c^{\tau},d^{\tau})w_{0}$$ for $c, d \in \cC$.
Therefore, the flip $\tau$ acts on the Dynkin diagram via conjugation
with the longest word $w_{0}$ of the Weyl group. Hence, if $\cT = (\cB_{+},\cB_{-},\dl_*)$ is a spherical twin building, then any ad\-ja\-cen\-cy-preserving involution $\tau$ that interchanges
$\cB_{+}$ and $\cB_{-}$ and maps some chamber onto an opposite chamber
is a flip if and only if the induced map $\widehat\tau$ on the building
$\cB=(\cC,\dl)$ satisfies
$\dl(c,d)=w_{0}\dl(c^{\widehat\tau},d^{\widehat\tau})w_{0}$ for all chambers
$c,d \in \cC$ where $w_{0}$ is the longest word in the Weyl group $W$.

\subsubsection{Flips and polarities} \label{flippol}

For a flip $\tau$ of a spherical twin building of type $A_n$ considered as the building geometry, i.e., the projective geometry, the induced map $\widehat\tau$ (see Section \ref{fg}) is an incidence-preserving involution that maps points onto hyperplanes such that for points $p$, $q$ one has an incidence between $p$ and $q^{\widehat\tau}$ if and only if $q$ and $p^{\widehat\tau}$ are incident. Hence $\widehat\tau$ is a polarity of the projective geometry, cf.\ \cite{Buekenhout/Shult:1974}. This means, by \cite{Buekenhout/Shult:1974}, \cite{Tits:1974}, \cite{Veldkamp:1959}, \cite{Veldkamp:1960}, that $\widehat{\tau}$ is induced by a pseudo-quadratic or an alternating form, if $n \geq 4$, see also \cite{Cohen:1995}. Therefore a flip is a natural generalisation of a polarity, and we are on our way towards a generalisation of Aschbacher's geometry for arbitrary twin buildings. 

I mention in passing that flipflop geometries coming from flips inducing non-degenerate symmetric bilinear forms have been studied in \cite{Altmann:2003}, \cite{Altmann/Gramlich:2006}, \cite{Gramlich/Maldeghem:2006}, \cite{Roberts:2005}. Although a flip inducing a non-degenerate alternating form cannot be a Phan involution, one can still study the geometry of chambers having minimal co-distance from their image under that flip. This yields the geometry on hyperbolic lines of a symplectic polar space, which has been studied in contexts different from Phan theory in \cite{Cuypers:1994}, \cite{Gramlich:2004c}, \cite{Hall:1988}, \cite{Hall:1989}. In \cite{Blok/Hoffman:2008} this geometry has finally been investigated from the point of view of Phan theory, yielding interesting presentations of symplectic groups.

\subsection{Phan's second theorem and the classical Phan-type theorem}

\subsubsection{Weak Phan systems of arbitrary spherical type} \label{weakphan2}

Let $\Delta$ be an irreducible spherical Coxeter diagram of rank at least three. A group $G$ admits a {\em weak Phan system of type $\Delta$ over $\mathbb{F}_{q^2}$}, if
$G$ contains subgroups $U_\alpha\cong \SL_2(q^2) \cong \SU_2(q^2)$, $\alpha \in \Delta$, and
$U_{\alpha,\beta}$, $\alpha, \beta \in \Delta$, so that the following hold:
\begin{description}
\item[{\bf (wP1)}] if $\node_\alpha\hspace{3em}\node_\beta$, then $[x,y]=1$ for all $x\in U_\alpha$ and
$y\in U_\beta$,
\item[{\bf (wP2)}] $U_{\alpha,\beta} \cong \left\{
\begin{array}{ll}
\mathrm{(P)SU}_3(q^2), & \mbox{in case $\node_\alpha\stroke{}\node_\beta$,} \\
\mathrm{(P)Sp}_{4}(q), & \mbox{in case $\node_\alpha\dstroke{}\node_\beta$;}
\end{array}\right.$ \\
moreover, $U_\alpha$ and $U_{\beta}$ form a standard pair (see below) in $U_{\alpha,\beta}$, and
\item[{\bf (wP3)}] $G = \gen{U_{\alpha,\beta} \mid \alpha, \beta \in \Delta}$.
\end{description}
For $U_{\alpha,\beta} \in \{ \SU_3(q^2), \Sp_4(q) \}$ define $G_{\alpha,\beta} := \SL_3(q^2), \Sp_4(q^2)$ accordingly.
A {\it standard pair} in $U_{\alpha,\beta}$ is a pair of subgroups isomorphic to
$\SU_2(q^2) \cong \SL_2(q)$ conjugate as a pair to the intersections $G_\alpha \cap U_{\alpha,\beta}$ and $G_\beta \cap U_{\alpha,\beta}$, where $G_\alpha$, $G_\beta$ form a pair of fundamental rank one subgroups of $G_{\alpha,\beta}$ (Section \ref{redundancy2}).  Standard pairs in central quotients are defined as the
images under the canonical homomorphism of standard pairs of the simply connected group. A concrete description of standard pairs of $\SU_3(q^2)$ can be found in Section \ref{weakphan}. For a concrete description of standard pairs of $\Sp_4(q)$ see \cite{Gramlich/Hoffman/Shpectorov:2003}, \cite{Gramlich/Horn/Nickel:2006}.

\subsubsection{Phan-type theorem of type $B_n$} \label{522}

The analogue of Aschbacher's geometry can be constructed from Example 2a (Section \ref{exclass}) as Example 1b (Section \ref{nondegspacerev}) has been deduced from Example 1a (Section \ref{exclass}). 

\medskip \noindent 
{\bf Example 2b.} Consider the situation of Example 2a, but with $\mathbb{F} = \mathbb{F}_{q^2}$, let $G \cong \Omega_{2n+1}(q^2)$, i.e., the commutator subgroup of $\mathrm{GO}_{2n+1}(q^2)$, and denote the form on $V$ by $\form$. Since the case of even $q$ will be covered in Section \ref{Cnphan} via the isomorphism $\Spin_{2n+1}(2^e) \cong \Sp_{2n}(2^e)$, it suffices to study the case of $q$ odd.  Let
$\{e_1,\ldots,e_n,f_1,\ldots,f_n,x\}$ be a hyperbolic basis of the
orthogonal space $V$, so that $(e_i,f_j)=\dl_{ij}$, while $x$ with $(x,x)=1$ is orthogonal to each basis vector except itself.  Consider the
semilinear transformation $\tau$ of $V$ which is the composition of the
linear transformation given by the Gram matrix of $\form$ with respect to the above basis and the
involutory field automorphism applied to the coordinates.  

It can be shown, cf.\ \cite{Bennett/Gramlich/Hoffman/Shpectorov:2007}, that $\tau$ produces a Phan involution of
$\cT$.  Furthermore, $\cC_\tau$ is geometric and $G_\tau\cong \SO_{2n}(q)$ (cf.\ \cite[Proposition 2.10]{Bennett/Gramlich/Hoffman/Shpectorov:2007})
acts flag-transitively on the corresponding flipflop geometry
$\G_\tau$.  The geometry $\G_\tau$ can be described as follows.  For
$u,v\in V$ let $((u,v))=(u,v^\tau)$.  Then $\fform$ is a non-degenerate hermitian form.  The
flipflop geometry $\G_\tau$ can be identified via $(x_+,x_-)\mapsto
x_+$ with the geometry $\G_{B_n}$ of all subspaces of $V$ which are totally
isotropic with respect to $\form$ and, at the same time, non-degenerate
with respect to $\fform$.

\medskip
In \cite{Bennett/Gramlich/Hoffman/Shpectorov:2007}, \cite{Gramlich/Horn/Nickel:2007} the simple connectedness of $\G_{B_n}$ is proved, leading to the following result. 

\begin{phantype} [Bennett, G., Hoffman, Shpectorov \cite{Bennett/Gramlich/Hoffman/Shpectorov:2007}, G., Horn, Nickel \cite{Gramlich/Horn/Nickel:2007}]
Let $q$ be an odd prime power, let $n \geq 3$, and let $G$ be a group admitting a weak Phan system of type $B_n$ over $\mathbb{F}_{q^2}$.
\begin{enumerate}
\item If $q\ge 5$, then $G$ is isomorphic to a quotient of $\Spin(2n+1,q)$.
\item For $n\ge 4$, let $G$ be a group
admitting a weak Phan system of type $B_n$ over $\mathbb{F}_9$. In
addition, assume that $\gen{U_{i-1},U_i,U_{i+1}}$ is isomorphic to
a central quotient of $\SU(4,9)$ (if $2\le i\le n-2$) or
$\Spin(7,3)$ (if $i=n-1$). Then $G$ is isomorphic to $\Spin(2n+1,3)$
or a central quotient thereof.
\end{enumerate}
\end{phantype}

\subsubsection{The group $\Omega(7,3)$}

In \cite{Gramlich/Horn/Nickel:2007} a group $H$ admitting a weak Phan system of type $B_3$ over $\F_{3^2}$ is constructed which is a $2187$-fold extension of $\Omega_7(3)$. 
To be precise, see \cite{Gramlich/Horn/Nickel:2007}, the group $H$ is isomorphic to a non-split extension of $\Omega(7,3)$ by $K:=(\mathbb{Z}/3\mathbb{Z})^7$, i.e. the sequence $1 \to K \to H \rightarrow \Omega(7,3) \to 1$ is exact and non-split. This extension of $\Omega_7(3)$ has been studied in \cite{Fischer:1971}, \cite{Kuesefoglu:1979}, \cite{Kuesefoglu:1980}. 
Altogether we can conclude that for $q=3$ the geometry $\G_{B_3}$ admits a 2187-fold covering, whence is not simply connected. It is shown in \cite{Gramlich/Horn/Nickel:2007} that this covering is universal.

\subsubsection{Phan-type theorem of type $C_n$} \label{Cnphan}

The geometry needed to prove the Phan-type theorem of type $C_n$ looks very much like the one of type $B_n$.

\medskip \noindent
{\bf Example 3b.} Consider the situation of Example 3a, but with the
field of definition of order $q^2$, so that $G \cong \Sp_{2n}(q^2)$.  Let
$\{e_1,\ldots,e_n,f_1,\ldots,f_n\}$ be a hyperbolic basis of the
symplectic space $V$, i.e., $(e_i,f_j)=\dl_{ij}$, where $\form$ is the alternating
form on $V$, and consider the
semilinear transformation $\tau$ of $V$ which is the composition of the
linear transformation given by the Gram matrix of $\form$ with respect to the above basis and the
involutory field automorphism applied to the coordinates.  

It can be shown, cf.\ \cite{Gramlich/Hoffman/Shpectorov:2003}, that $\tau$ produces a Phan involution of
$\cT$.  Furthermore, $\cC_\tau$ is geometric and $G_\tau\cong \Sp_{2n}(q)$
acts flag-transitively on the corresponding flipflop geometry
$\G_\tau$.  By \cite{Gramlich/Hoffman/Shpectorov:2003} the geometry $\G_\tau$ has the following alternative description.  For
$u,v\in V$ let $((u,v))=(u,v^\tau)$, so that $\fform$ is a non-degenerate hermitian form.  The
flipflop geometry $\G_\tau$ can be identified via $(x_+,x_-)\mapsto
x_+$ with the geometry $\G_{C_n}$ of all subspaces of $V$ which are totally
isotropic with respect to $\form$ and, at the same time, non-degenerate
with respect to $\fform$.

\medskip
By \cite{Gramlich/Hoffman/Shpectorov:2003} (with the missing cases dealt with in \cite{Gramlich/Horn/Nickel:2006}, \cite{Horn:2005}) the geometry $\G_\tau$ is almost always simply connected, resulting in the following Phan-type theorem. 

\begin{phantype}[G., Hoffman, Shpectorov \cite{Gramlich/Hoffman/Shpectorov:2003}, G., Horn, Nickel \cite{Gramlich/Horn/Nickel:2006}, Horn \cite{Horn:2005}] \label{main3} \label{main4}
Let $q$ be a prime power, let $n \geq 3$, and let $G$ be a group admitting a weak Phan system of type $C_n$ over $\mathbb{F}_{q^2}$.
\begin{enumerate}
\item If $q \geq 3$, then $G$ is isomorphic to a central quotient of $\Sp_{2n}(q)$.
\item If $q=2$ and $n\geq4$ and if
\begin{enumerate}
\item for any triple $i,j,k$ of nodes of the Dynkin diagram $C_n$ that form a subdiagram
\[ \node^i\arc\node^j\arc\node^k \]
of type $A_3$, the subgroup $\langle U_{i,j},U_{j,k}\rangle$ is isomorphic to a central quotient of $\SU_4(2^2)$;
\item for any triple $i,j,k$ of nodes of the Dynkin diagram $C_n$ that form a subdiagram
\[ \node^i\arc\node^j\dstroke{<}\node^k \]
of type $C_3$, the subgroup $\langle U_{i,j},U_{j,k}\rangle$ is isomorphic to a central quotient of $\Sp_6(2)$;
\item
  \begin{enumerate}
  \item for any triple $i,j,k$ of nodes of the Dynkin diagram $C_n$ that form a subdiagram
\[ \node^i\qquad\quad\node^j\arc\node^k \]
of type $A_1\oplus A_2$, the groups $U_i$ and $U_{j,k}$ commute elementwise; and

  \item for any quadruple 
   of nodes of the Dynkin diagram $C_n$ that form a subdiagram
\[ \node^i\arc\node^j\qquad\quad\node^k\arc\node^l \]
of type $A_2\oplus A_2$, the groups $U_{i,j}$ and $U_{k,l}$ commute elementwise; and

  \item for any triple $i,j,k$ of nodes of the Dynkin diagram $C_n$ that form a subdiagram
\[ \node^i\qquad\quad\node^j\dstroke{<}\node^k \]
of type $A_1\oplus C_2$, the groups $U_i$ and $U_{j,k}$ commute elementwise; and

  \item for any quadruple 
   of nodes of the Dynkin diagram $C_n$ that form a subdiagram
\[ \node^i\arc\node^j\qquad\quad\node^k\dstroke{<}\node^l \]
of type $A_2\oplus C_2$, the groups $U_{i,j}$ and $U_{k,l}$ commute elementwise;

  \end{enumerate}
\end{enumerate}
then $G$ is isomorphic to a central quotient of $\Sp_{2n}(2)$.
\end{enumerate}
\end{phantype}

Unlike the cases $A_3$ and $B_3$, the geometry $\G_{C_3}$ is simply connected even for $q=3$, cf.\ \cite{Gramlich/Horn/Nickel:2006}, \cite{Horn:2005}.

\subsubsection{Phan's second theorem and Phan-type theorem of type $D_n$} \label{phan2}

The last series $D_n$ of finite classical groups, the even-dimensional orthogonal groups, again belongs to a simply laced diagram, which has already been treated in \cite{Phan:1977a}. The arguments in \cite{Phan:1977a} are based on the construction of a presentation that identifies the target group as an orthogonal group via \cite{Wong:1974}. Here is the main result of \cite{Phan:1977a} concerning $D_n$.

\begin{phan}[\cite{Phan:1977a}] \label{Phan2}
Let $q \geq 5$ be odd and let $n \geq 4$. If $G$ admits a Phan system of type $D_n$ over $\Fqsq$,
then $G$ is isomorphic to a factor group of $\Spin_{2n}^+(q^2)$, if $n$ is even, and isomorphic to a factor group of $\Spin_{2n}^-(q^2)$, if $n$ is odd.
\end{phan}

This result has been revised in \cite{Gramlich/Hoffman/Nickel/Shpectorov:2005} using the following geometry.

\medskip
\noindent
{\bf Example 4b.} Consider the situation as in Example 4a, but over the field $\Fqsq$, and let $G = \Omega^+_{2n}(q^2)$. For sake of simplicity of the exposition we assume here that $q$ is odd, although in \cite{Gramlich/Hoffman/Nickel/Shpectorov:2005} also the case of even characteristic is dealt with. The Phan involution $\tau$ can again be defined as the
composition of the linear transformation given by the Gram matrix of the bilinear form $\form$ with respect to a hyperbolic basis and coordinate-wise application of the involutory field automorphism.  This $\tau$ produces a
flipflop geometry on which $G_\tau \cong \Omega^\pm_{2n}(q)$ acts flag-transitively, cf.\ \cite[Proposition 3.10]{Gramlich/Hoffman/Nickel/Shpectorov:2005}. The geometry $\G_\tau$ can be described as follows.  For
$u,v\in V$ let $((u,v))=(u,v^\tau)$, where $\form$ is the orthogonal
form on $V$, so that $\fform$ is a non-degenerate hermitian form.  The
flipflop geometry $\G_\tau$ can be identified via $(x_+,x_-)\mapsto
x_+$ with the geometry $\G_{D_n}$ of all subspaces of $V$ which are totally
isotropic with respect to $\form$ and, at the same time, non-degenerate
with respect to $\fform$. See \cite{Gramlich/Hoffman/Nickel/Shpectorov:2005} for more details and a description of the geometry for even $q$.

\begin{phantype}[G.\ Hoffman, Nickel, Shpectorov \cite{Gramlich/Hoffman/Nickel/Shpectorov:2005}] \label{main5} \label{main6}
Let $q$ be a prime power, let $n \geq 3$, and let $G$ be a group admitting a weak Phan system of type $D_n$ over $\F_{q^2}$.
\begin{enumerate}
\item If $q \geq 4$, then $G$ is isomorphic to a central quotient of
\begin{itemize}
\item $\Spin^+_{2n}(q)$, if $n$ even; and
\item $\Spin^-_{2n}(q)$, if $n$ odd.
\end{itemize}
\item If $q = 2$, $3$ and $n \geq 4$ and if, furthermore,
\begin{enumerate}
\item for any triple $i$, $j$, $k$ of nodes of the Dynkin diagram $D_n$ that form a subdiagram $$\node^{i}\stroke{}\node^{j}\stroke{}\node^{k}$$ of type $A_3$, the subgroup 
$\gen{U_{i,j},U_{j,k}}$ is isomorphic to a central quotient of 
$\SU_4(q^2)$;
\item in case $q=2$
\begin{enumerate}
\item for any triple $i$, $j$, $k$ of nodes of the Dynkin diagram $D_n$ that form a subdiagram $$\node^{i} \quad \quad \node^{j}\stroke{}\node^{k}$$ of type $A_1 \oplus A_2$ the groups $U_i$ and $U_{j,k}$ commute 
elementwise; and 
\item for any quadruple $i$, $j$, $k$, $l$ of nodes of the Dynkin diagram $\Delta$ that form a subdiagram $$\node^{i}\stroke{}\node^{j} \quad \quad \node^{k}\stroke{}\node^{l}$$ of type $A_2 \oplus A_2$
the groups
$U_{i,j}$ and $U_{k,l}$ commute elementwise;
\end{enumerate}
\end{enumerate}
then $G$ is isomorphic to a central quotient of
\begin{itemize}
\item $\Spin^+_{2n}(q)$, if $n$ even; and
\item $\Spin^-_{2n}(q)$, if $n$ odd.
\end{itemize}
\end{enumerate}
\end{phantype}

\subsection{The Devillers-M\"uhlherr filtration} \label{DM}

\subsubsection{Filtrations of chamber systems}

Lacking concrete easy models of the flipflop geometries of exceptional type some extra theory is necessary in order to be able to extend the Phan-type theorems to exceptional type groups.

Such a theory has been developed in \cite{Devillers/Muehlherr:2007}:
A \textit{filtration} of a chamber system ${\cal C} = (C,(\sim_i)_{i \in I})$ over $I$ is a family 
${\cal F} = (C_n)_{n \in \mathbb{N}}$
of subsets of $C$ such that
\begin{enumerate}
\item $C_n \subset C_{n+1}$ for all $n \in \mathbb{N}$,
\item $\bigcup_{n \in \mathbb{N}} C_n = C$, and
\item for each $n >0$ with $C_{n-1}\neq \emptyset$ there exists $i \in I$
such that for each chamber $c \in C_n$ there exists a chamber
$c' \in C_{n-1}$ which is $i$-adjacent to $c$.
\end{enumerate}
A filtration ${\cal F} = (C_n)_{n \in \mathbb{N}}$ is called
\textit{residual}, if for each $\emptyset \neq J \subset I$ and each
$c \in C$ the family $(C_n \cap R_J(c))_{n \in \mathbb{N}}$
is a filtration of the chamber system 
${\cal R}_J(c) = (R_J(c), (\sim_j)_{j \in J})$.
For $x \in C$ define $|x| := \min\{ \lambda \in \mathbb{N} 
\mid x \in C_{\lambda} \}$ and, for $X \subseteq C$,
define $|X| := \min\{ |x| \mid x \in X \}$ and $\aff(X) := 
\{ x \in X \mid |x| = |X| \}$.
Note that $\aff(C) = C_m$
where $m = \min \{ n \in \mathbb{N} \mid C_n \neq \emptyset \}$.

\subsubsection{Filtering flipflop systems inside buildings} \label{ffsib}

Let $\mathcal{T} = ({\cal B}_+,{\cal B}_-,\delta_*)$ be a twin building of type $(W,S)$ with a flip $\tau$ satisfying axiom (F4)' from Section \ref{flipphan}. Then there exists a filtration $${\cal F}_\tau = (C_n)_{n \in \mathbb{N}}$$ of the building ${\cal B}_+$ so that $C_0$ equals the set of chambers of the flipflop system $\mathcal{C}_\tau$ defined as follows.

For a residue $R$ of $\mathcal{B}_+$ put $l_*(\tau,R) := \min \{ l(\delta_*(x,\tau(x))) \mid x \in R \}$
and $A_{\tau}(R) := \{ x \in R \mid  l(\delta_*(x,\tau(x))) = l_*(\tau,R)\}$, where $l$ denotes the length function of the group $W$ with respect to the generating set $S$. Since $S$ is finite, there exists an injective map $\Inv(W) \to \mathbb{N} : x \mapsto |x|$ from the involutions of $W$ to the non-negative integers with $|1_W| = 0$ such that $l(x) < l(y)$ implies
$|x| < |y|$. Defining $$C_n := \{ c \in C_+ \mid |\delta_*(c,\tau(c))| \leq n \},$$ the family ${\cal F}_{\tau} = (C_n)_{n \in \mathbb{N}}$ is a residual filtration
of ${\cal C}({\cal B}_+)$ by \cite{Devillers/Muehlherr:2007}.

\subsubsection{A criterion for simple connectedness of a flipflop system} \label{sccrit}

The setup from Section \ref{ffsib} and the simple $2$-connectedness of buildings (cf.\ the Solomon-Tits Theorem in Section \ref{buildings}) yield the following criterion of simple connectedness of flipflop systems established in \cite{Devillers/Muehlherr:2007}.

If $\tau$ is a flip satisfying axiom (F4)' of a three-spherical twin building ${\cal T} = (\mathcal{B}_+,\mathcal{B}_-,\delta_*)$ of finite rank (i.e., a twin building of finite rank whose residues of rank three are spherical) such that
\begin{enumerate}
\item  the chamber system  $(A_{\tau}(R), (\sim_t)_{t \in J})$ is connected for each $J$-residue $R$ of rank two, and
\item the chamber system $(A_{\tau}(R), (\sim_t)_{t \in J})$ is simply 2-connected for each $J$-residue $R$ of rank three,
\end{enumerate}
then the flipflop system $\mathcal{C}_\tau$ is simply 2-connected in the sense of \ref{2-connected}.

\subsection{Wedges of spheres and the Abels--Abramenko filtration} \label{soltitsa}

\subsubsection{Generalised flipflop geometries of type $A_n$} \label{descriptionan}

In view of Section \ref{sccrit} it remains to study the chamber systems $$(A_{\tau}(R), (\sim_t)_{t \in J})$$ for residues $R$ of rank two and three in order to prove the simple connectedness of the exceptional flipflop geometries.
In case the diagram of the twin building $\mathcal{T}$ is simply laced, these chambers systems can be described by so-called generalised flipflop geometries of type $A_n$, defined in this section, cf.\ \cite{Blok/Hoffman:2009}, \cite{Devillers/Gramlich/Muehlherr}, \cite{Gramlich/Witzel}.

Two subspaces $A$ and $B$ of a vector space $V$ are {\em opposite} when $V=A\oplus B$. A subspace $A$ is {\em transversal} or {\em in general position} to a flag $F$, i.e., a chain of incident subspaces of $V$, if for any subspace $B$ of $F$ we have $A\cap B=\{0\}$ or $V=A+B$.  
In other words, $A$ is transversal to $F$, in symbols $A\pitchfork_V F$, if and only if there is a subspace $C$ of $V$ incident with $F$ such that $A$ and $C$ are opposite.

For a field $\mathbb{F}$ with an involution $\sigma$ and an $(n+1)$-dimensional $\mathbb{F}$-vector space $V$ containing a flag $F$ equal to $0 = V_0 \lneq V_1 \lneq \cdots \lneq V_t \lneq V_{t+1} = V$ of subspaces of $V$ endowed with $\sigma$-hermitian forms $\omega_{i} : V_{i+1} \times V_{i+1} \to \mathbb{F}$, $0 \leq i \leq t$, satisfying $\Rad\left(\omega_i\right) = V_i$, the {\em generalised flipflop geometry of type $A_n$} (modelled in $V$ with respect to the flag $F$ and the forms $\omega_i$) consists of all proper non-trivial vector subspaces $U$ of $V$ transversal to $F$ with $U \cap V_{k_U+1}$ non-degenerate with respect to $\omega_{k_U}$ where $k_U = \min \{ i \in \{0, \ldots, t \} \mid U \cap V_{i+1} \not= \{0\} \}$. 

In the simply laced three-spherical case over $\mathbb{F} = \mathbb{F}_{q^2}$ a geometry arising from a chamber system $(A_{\tau}(R), (\sim_t)_{t \in J})$, $|J| \in \{ 2, 3 \}$ (defined in Section \ref{ffsib}), is isomorphic to a generalised flipflop geometry for $n \in \{ 2, 3 \}$ by \cite{Blok/Hoffman:2009}, \cite[Proposition 6.6]{Gramlich/Witzel}.

For $t = 0$ and $\mathbb{F} = \mathbb{F}_{q^2}$, the generalised flipflop geometry on $V$ equals Aschbacher's geometry on $V$, i.e., the flipflop geometry of type $A_n$ over $\mathbb{F}_{q^2}$, cf.\ \cite{Bennett/Gramlich/Hoffman/Shpectorov:2003}, \cite{Bennett/Shpectorov:2004} and Section \ref{nondegspacerev}.

For $t = n$, the generalised flipflop geometry on $V$ equals the geometry opposite the chamber $F$. This follows from the fact that each ${\omega_i}$ has rank one with radical $V_i$. Therefore any vector $v \in V_{i+1} \backslash V_i$ is non-degenerate with respect to $\omega_i$, so that any subspace $U$ of $V$ with $U \oplus V_i = V$ intersects $V_{i+1}$ in a non-degenerate (with respect to $\omega_i$) one-dimensional subspace. 

\subsubsection{A Solomon--Tits-type Theorem} \label{wedgesspheres}

It turns out that generalised flipflop geometries of type $A_n$ are not only a useful tool in order to prove Phan-type theorems for groups with simply laced diagrams, but are also interesting in their own right. Indeed, via the Abels--Abramenko filtration \cite{Abels/Abramenko:1993} it can be shown that a generalised flipflop geometry $\G$ of type $A_n$ is homotopy equivalent to a wedge of $(n-1)$-spheres provided the field $\mathbb{F}$ contains sufficiently many elements. 

In order to describe this filtration let $p$ be a one-dimensional subspace of $V$ which is non-degenerate with respect to the hermitian form $\omega_t$ and define $Y_0:=\{W\in \G \mid \gen{p,W} \in \G\}$ and $Y_i:=Y_{i-1}\cup\{W\in \G \mid \dim W=n+1-i\}$ for $1 \leq i \leq n$. The strategy from \cite{Abels/Abramenko:1993} can be transferred literally to obtain the following Solomon--Tits-type result, which (similar to what is surveyed in \cite{Humphreys:1987}) gives rise to a representation of $\SU_{n+1}(\mathbb{F}_{q^2})$ on the integral homology group $H_{r-1}(\G)$ tensored with $\mathbb{Q}$, which may be an interesting object to study.

\begin{solomontits}[Devillers, G., M\"uhlherr \cite{Devillers/Gramlich/Muehlherr}] \label{allrad}
Let $V$ be an $(n+1)$-dimensional vector space over a field $\mathbb{F}$ with an involution, let $(\G_j)_{1 \leq j \leq m}$ be a finite family of generalised flipflop geometries of type $A_n$ modelled in $V$, and let $\G = \bigcap_j \G_j$. In case $\mathbb{F} = \mathbb{F}_{q^2}$ assume  $2^{n-1}(q+1)m < q^2$. Then $|\G|$ is homotopy equivalent to a wedge of $(n-1)$-spheres. 
\end{solomontits}

Notice in passing that this result once again proves simple connectedness of Aschbacher's geometry, at least for large fields. 

Moreover, this result can be used to deduce finiteness properties of the group $\SU_{n+1}(\mathbb{F}_{q^2}[t,t^{-1}],\theta)$ in the spirit of \cite{Abramenko:1996}, \cite[Chapter 13]{Abramenko/Brown:2008}, where $\theta$ is the involution of $\SL_{n+1}(\mathbb{F}_{q^2}[t,t^{-1}])$ which acts as the Chevalley involution on $\SL_{n+1}$, as the Frobenius involution on $\mathbb{F}_{q^2}$, and interchanges $t$ and $t^{-1}$. In fact, this group $\SU_{n+1}(\mathbb{F}_{q^2}[t,t^{-1}],\theta)$ is a lattice in $\SL_{n+1}(\mathbb{F}_{q^2}((t)))$ and in $\SL_{n+1}(\mathbb{F}_{q^2}((t^{-1})))$, cf.\ \cite{Gramlich/Muehlherr:2008}, whence an arithmetic group by \cite[Chapter IX]{Margulis:1991}. See \cite{Blok/Hoffman:2009} for a concrete description of the group $\SU_{n+1}(\mathbb{F}_{q^2}[t,t^{-1}],\theta)$ and related groups.

\subsection{Phan's third theorem and the Phan-type theorem of type $E_n$}

\subsubsection{Phan's Theorem}

The article \cite{Phan:1977a} also contains a theorem concerning the diagrams $E_6$, $E_7$, and $E_8$. Phan's Theorem \ref{Phan2} (Section \ref{phan2}) plus \cite{Phan:1977} are used in order to construct a system of subgroups satisfying the hypotheses of the Curtis--Tits Theorem Version \ref{phanstyle}, which then is invoked.

\begin{phan}[Phan \cite{Phan:1977a}] \label{Phan3}
Let $q \geq 5$ be odd. If $G$ admits a Phan system of type $E_6$, $E_7$, or $E_8$ over $\Fqsq$,
then $G$ is isomorphic to a factor group of the universal Chevalley group $^2E_6(q^2)$, $E_7(q)$, or $E_8(q)$, respectively.
\end{phan}

\subsubsection{Exploiting the filtrations} \label{exploiting}

By the Solomon-Tits-type Theorem \ref{allrad} (Section \ref{wedgesspheres}) a generalised flipflop geometry of type $A_3$ over $\mathbb{F}_{q^2}$ is simply connected, provided $2^2(q+1) < q^2$, which is the case for $q \geq 5$, while a generalised flipflop geometry of type $A_2$ over $\mathbb{F}_{q^2}$ is connected, if $2(q+1) < q^2$, which is the case for $q \geq 3$. Together with the criterion for simple connectedness of a flipflop system from \cite{Devillers/Muehlherr:2007} (Section \ref{sccrit}) this implies that the flipflop geometries of type $E_6$, $E_7$, $E_8$ over $\mathbb{F}_{q^2}$ are simply connected provided $q \geq 5$. 

For completeness I should point out here that the chamber systems $$(A_\tau(R),(\sim_t)_{t \in J})$$ for residues $J$ of type $A_1 \oplus A_1$, $A_1 \oplus A_2$, $A_1 \oplus A_1 \oplus A_1$ are automatically (simply) connected by the following standard argument. Assuming that $\G = \G_1 \oplus \G_2$ with $\G_1$ connected of rank 
at least two and $\G_2$ non-empty, the geometry $\G$ is simply connected.
Indeed, the geometry $\G$ is certainly connected, and choosing a base point $x \in \G_1$ one can prove that any cycle originating at $x$ is homotopic to a cycle 
fully contained in $\G_1$. Such a cycle then is null homotopic because it forms a cone together with any element $z \in \G_2$.  

\subsubsection{The Phan-type theorem of type $E_n$}

Alternatively | and this had already been done by Hoffman, M\"uhlherr, Shpectorov and the author roughly one year before the Solomon-Tits-type Theorem \ref{allrad} was proved | one can directly compute the fundamental group for generalised flipflop geometries of type $A_3$. It turns out that via direct computation it is possible to show that the fundamental groups are trivial for $q \geq 4$. Together with the classification of amalgams \cite{Bennett/Shpectorov:2004}, \cite{Dunlap:2005} and the criterion for simple connectedness of a flipflop system in \cite{Devillers/Muehlherr:2007} (Section \ref{sccrit}) this implies the following Phan-type theorem.

\begin{phantype}[G., Hoffman, M\"uhlherr, Shpectorov 2005]
Let $q\ge 4$ be a prime power and let $G$ be
a group containing a weak Phan system of type $E_6$, $E_7$, or $E_8$ over $\Fqsq$.
Then $G$ is isomorphic to a group of type ${}^2E_6(q^2)$, $E_7(q)$, or $E_8(q)$.
\end{phantype}

\subsection{The Abramenko filtration and the Phan-type theorem of type $F_4$}

\subsubsection{Generalised flipflop geometries of type $B_n$ and $C_n$} \label{descriptioncn}

The criterion from Section \ref{sccrit} allows three-spherical diagrams. In view of the method of proof of the Phan-type theorem of type $E_n$ via generalised flipflop geometries of type $A_3$ it is natural to ask for the definition of generalised flipflop geometries of type $B_3$ and $C_3$.

Let $V$ be a vector space over a field $\mathbb{F}$ with an involution $\sigma$. For $U,W \le V$, we say that $U$ is \emph{transversal} to $W$ and write $U \pitchfork W$, if $U \cap W = 0$ or $\gen{U,W} = V$. Note that $U \pitchfork W$ if and only if $\dim(U \cap W) = \max\{0,\dim U +\dim W - \dim V\}$.
For a flag $F=(0=V_0 \le \ldots \le V_k=V)$ and a subspace $U \le V$ we say that $U$ is \emph{transversal} to $F$ and write $U \pitchfork F$, if $U \pitchfork V_i$ for $0 \le i \le k$. This is the case if and only if $\gen{U,V_{k_U}} = V$ where $k_U = \min \{i \mid U \cap V_j \ne \{0\}\}$.

Given a flag $F=(0=V_0 \le \ldots \le V_k=V)$ we call a family $(\omega_i)_{1 \le i \le k}$ of $\sigma$-hermitian forms $\omega_i \colon V_i \times V_i \to \mathbb{F}$ \emph{compatible with $F$} if $\mathrm{Rad}(\omega_i) = V_{i-1}$.

Let $F$ be as above and let $\omega=(\omega_i)_i$ be a family of compatible $\sigma$-hermitian forms. For $U \le V$ we say that $U$ is \emph{transversal} to $(F,\omega)$, if $U$ is transversal to $F$ and $U \cap V_{k_U}$ is $\omega_{k_U}$-non-degenerate. In this case we write $U \pitchfork (F,\omega)$.

Let $\Delta$ be the building geometry of type $B_n(\mathbb{F})$ or $C_n(\mathbb{F})$ embedded in a $\mathbb{F}$-vector space $V$ of dimension $2n+1$, resp.\ $2n$, and let $e_1$, \ldots, $e_n$, $f_1$, \ldots, $f_n$, $x$ be a standard hyperbolic basis of $V$ (where the vector $x$, of course, only occurs in case $B_n$).
Let $F=(0=V_0 \le \ldots \le V_k=V)$ be a flag satisfying $F^\perp = F$. Let $\omega$ be a family of $\sigma$-hermitian forms compatible with $F$ and assume that there is an $\omega_k$-non-isotropic vector that is $(\cdot,\cdot)$-isotropic. The {\em generalised flipflop geometry} of type $\Delta$ over $\mathbb{F}$ defined by $(F,\omega)$ consists of all subspaces $U$ of $V$ that are totally $(\cdot,\cdot)$-isotropic and transversal to $(F,\omega)$.

A closer look reveals that half of the forms $\omega_i$ actually do not play any role, because a totally isotropic subspace $U$ that is transversal to $F$ cannot meet any of the $V_i$ with $\dim V_i \le n$. However, taking this into account would not simplify anything, but would instead make the definition of a generalised Phan geometry even more cumbersome.

\subsubsection{Another Solomon-Tits-type theorem and the Phan-type theorem of type $F_4$} \label{ansoltits}\label{F4}

The concept of transversality introducted in Section \ref{descriptioncn} makes the Abramenko filtration from \cite{Abramenko:1996} accessible. This filtration has been used in \cite{Gramlich/Witzel} in order to prove the following theorem.

\begin{solomontits}[G., Witzel \cite{Gramlich/Witzel}]
Let $\mathbb{F}$ be a field with an involution $\sigma$, let $(\G_j)_{1 \leq j \leq m}$ be a finite family of generalised flipflop geometries of type $B_n$ or $C_n$ embedded in some $(2n+1)$- or $2n$-dimensional $\mathbb{F}$-vector space $V$, and let $\G = \bigcap_j \G_j$. In case $\mathbb{F} = \mathbb{F}_{q^2}$ assume  $4^{n-1}(q+1)m < q^2$. Then $|\G|$ is homotopy equivalent to a wedge of $(n-1)$-spheres. 
\end{solomontits}

Similar to the case $A_3$ Hoffman, M\"uhlherr, Shpectorov and the author proved by direct computation that generalised flipflop geometries of type $B_3$ or $C_3$ are simply connected provided the underlying field contains at least $13$ elements. Again using the simple connectedness criterion from \cite{Devillers/Muehlherr:2007} (Section \ref{sccrit}), the final Phan-type theorem follows. Note that the generalised flipflop geometries are the correct objects in order to describe the chamber system $(A_\tau(R),(\sim_t)_{t \in J})$ by \cite[Proposition 6.6]{Gramlich/Witzel}. 

\begin{phantype}[G., Hoffman, M\"uhlherr, Shpectorov 2007]
Let $q\ge 13$ be a prime power and let $G$ be
a group containing a weak Phan system of type $F_4$ over $\Fqsq$.
Then $G$ is isomorphic to a group of type $F_4(q)$.
\end{phantype}

It remains to study the cases of small $q$.  

\section{Statement of the Phan-type theorem over finite fields} \label{complete}

We have reached one of the main purposes of this survey, the statement of the Phan-type theorem over finite fields. From Section \ref{fg} we know for which groups of Lie type the Phan-type theorem can make a statement, namely ${}^2A_n$, $B_n$, $C_n$, $D_{2n}$, ${}^2D_{2n+1}$, ${}^2E_6$, $E_7$, $E_8$, $F_4$. 

\medskip
\noindent {\bf The Phan-type theorem for finite fields.} {\em Let $q \geq 3$, let $\Delta$ be a spherical Dynkin diagram of rank at least three, and let $G$ be a group with a weak Phan system of type $\Delta$ over $\mathbb{F}_{q^2}$. Then $G$ is isomorphic to a quotient of
\begin{itemize}
\item $\SU_{n+1}(q^2)$, if $\Delta = A_n$ and $q \geq 4$ \\ {\rm (Bennett, Shpectorov
\cite{Bennett/Shpectorov:2004}, Phan \cite{Phan:1977})};
\item $\Spin_{2n+1}(q)$, if $\Delta = B_n$ and $q \geq 4$ \\ {\rm (Bennett, G., Hoffman, Shpectorov \cite{Bennett/Gramlich/Hoffman/Shpectorov:2007}, G., Horn, Nickel \cite{Gramlich/Horn/Nickel:2007})};
\item $\Sp_{2n}(q)$, if $\Delta = C_n$ \\ {\rm (G., Hoffman, Shpectorov \cite{Gramlich/Hoffman/Shpectorov:2003}, G., Horn, Nickel \cite{Gramlich/Horn/Nickel:2006}, Horn \cite{Horn:2005})};
\item $\Spin_{2n}^\pm$, if $\Delta = D_n$ and $q \geq 4$, of plus type if $n$ even, of minus type if $n$ odd \\ {\rm (G., Hoffman, Nickel, Shpectorov
\cite{Gramlich/Hoffman/Nickel/Shpectorov:2005}, Phan \cite{Phan:1977a})};
\item the universal Steinberg--Chevalley group of type $^{2}E_6(q^2)$, if $\Delta = E_6$ and $q \geq 4$ \\ {\rm
(Devillers, G., M\"uhlherr \cite{Devillers/Gramlich/Muehlherr}, G., Hoffman, M\"uhlherr, Shpectorov
2005, Phan \cite{Phan:1977a})};
\item the universal Steinberg--Chevalley group of type $E_7(q)$, if $\Delta = E_7$ and $q \geq 4$ \\ {\rm
(Devillers, G., M\"uhlherr \cite{Devillers/Gramlich/Muehlherr}, G., Hoffman, M\"uhlherr, Shpectorov
2005, Phan \cite{Phan:1977a})};
\item the universal Steinberg--Chevalley group of type $E_8(q)$, if $\Delta = E_8$ and $q \geq 4$ \\ {\rm
(Devillers, G., M\"uhlherr \cite{Devillers/Gramlich/Muehlherr}, G., Hoffman, M\"uhlherr, Shpectorov
2005, Phan \cite{Phan:1977a})};
\item the universal Steinberg--Chevalley group of type $F_4(q)$, if $\Delta = F_4$ and $q \geq 13$ \\ {\rm (G., Hoffman, M\"uhlherr, Shpectorov 2007, G., Witzel \cite{Gramlich/Witzel})}.
\end{itemize}
}

\section{Curtis--Tits theory, Phan theory, and the revision of the classification of the finite simple groups} \label{app}

In this section I briefly mention by way of example how to prove a local recognition result for Chevalley groups of simply laced type. I currently do not know how to deal with the non-simply laced case.

Following \cite{Aschbacher:1977} (Sections \ref{phanasch}, \ref{redundancy2}) a {\em fundamental (rank one) subgroup} of a (twisted) Chevalley group $G$ is a group generated by two root subgroups $X_\alpha$, $X_{-\alpha}$, respectively the subgroup of fixed points of $\gen{X_\alpha,X_{-\alpha}}$ with respect to an involution of $G$ interchanging $X_\alpha$ and $X_{-\alpha}$.    

In the revision of the classification of the finite simple groups \cite{Gorenstein/Lyons/Solomon:1994}, \cite{Gorenstein/Lyons/Solomon:1995}, \cite{Gorenstein/Lyons/Solomon:1998}, \cite{Gorenstein/Lyons/Solomon:1999}, \cite{Gorenstein/Lyons/Solomon:2002}, \cite{Gorenstein/Lyons/Solomon:2006} one is interested in proving local recognition results of the following type.

\begin{locrec}[Altmann, G.\ 2007] \label{twisted}
Let $q$ be an odd prime power and let $G$ be a group containing an involution $x$ and a subgroup $K \unlhd C_G(x)$ such that 
\begin{enumerate}
\item $K \cong \left\{ \begin{array}{lc}
\SL_6(q) & (\mathrm{CT}) \\ \SU_{6}(q^2) & (\mathrm{P}) \end{array} \right.$;
\item $C_G(K)$ contains a subgroup $X \cong \SL_2(q) \cong \SU_2(q^2)$ with $\langle x \rangle = Z(X)$;
\item there exists an involution $g \in G$ such that $Y := gXg$ is contained in $K$; 
\item if $V$ is a natural module for $K$, then the commutator $[Y,V] = \{ yv - v \in V \mid y \in Y, v \in V \}$ is a subspace of $V$ of  $\left. \begin{array}{lc}
\mathbb{F}_{q} & (\mathrm{CT}) \\ \mathbb{F}_{q^2} & (\mathrm{P}) \end{array} \right\}$-dimension two;
\item $G = \gen{K,gKg}$; moreover, there exists $ z \in K \cap gKg$ which is a $gKg$-conjugate of $x$ and a $K$-conjugate of $gxg$.
\end{enumerate}
Then (up to isomorphism) 
\begin{gather}
G/Z(G)
\cong  {\rm PSL}_{8}(q) \mbox{ or } G/Z(G) \cong E_6(q), \tag{{\rm CT}} \\
G/Z(G)
\cong {\rm PSU}_{8}(q^2) \mbox{ or } G/Z(G) \cong {}^2E_6(q^2). \tag{{\rm P}}
\end{gather}
\end{locrec}

Using ideas developed in \cite{Cohen/Cuypers/Gramlich:2005}, \cite{Gramlich:2002}, \cite{Gramlich:2004b} the above theorem is implied by a graph-theoretical local recognition theorem. From the hypotheses of the theorem one constructs a connected locally line-hyperline graph (cf.\ \cite{Gramlich:2004b}), resp.\ a connected locally unitary line graph (cf.\ \cite{Altmann/Gramlich2}) with $G$ as a group of automorphisms and an induced subgraph $\Sigma$ isomorphic to the commuting reflection graphs $\bW(A_7)$ or $\bW(E_6)$ (see \cite{Gramlich/Hall/Straub}). This information then implies the existence of a Curtis--Tits, resp.\ Phan amalgam inside $G$ from which the theorem follows.

\addcontentsline{toc}{section}{Bibliography}

{\footnotesize
\bibliographystyle{plain}
\bibliography{habil}

\begin{thebibliography}{100}

\bibitem{Abels/Abramenko:1993}
Herbert Abels and Peter Abramenko.
\newblock On the homotopy type of subcomplexes of {Tits} buildings.
\newblock {\em Adv.\ Math.}, 101:78--86, 1993.

\bibitem{Abramenko:1996}
Peter Abramenko.
\newblock {\em Twin buildings and applications to {$S$}-arithmetic groups},
  volume 1641 of {\em Lecture Notes in Mathematics}.
\newblock Springer, Berlin, 1996.

\bibitem{Abramenko/Brown:2008}
Peter Abramenko and Kenneth~S.\ Brown.
\newblock {\em Buildings -- Theory and Applications}, volume 248 of {\em
  Graduate Texts in Mathematics}.
\newblock Springer, Berlin, 2008.

\bibitem{Abramenko/Muehlherr:1997}
Peter Abramenko and Bernhard M{\"u}hlherr.
\newblock Pr{\'e}sentation de certaines {$BN$}-paires jumel{\'e}es comme sommes
  amalgam{\'e}es.
\newblock {\em C.\ R.\ Acd.\ Sci.\ Paris S{\'e}r.\ I Math.}, 325:701--706,
  1997.

\bibitem{Abramenko/Ronan:1998}
Peter Abramenko and Mark~A.\ Ronan.
\newblock A characterization of twin buildings by twin apartments.
\newblock {\em Geom.\ Dedicata}, 73:1--9, 1998.

\bibitem{Abramenko/Maldeghem:1999}
Peter Abramenko and Hendrik Van~Maldeghem.
\newblock Connectedness of opposite-flag geometries in {Moufang} polygons.
\newblock {\em European J.\ Combin.}, 20:461--468, 1999.

\bibitem{Abramenko/Maldeghem:2000}
Peter Abramenko and Hendrik Van~Maldeghem.
\newblock On opposition in spherical buildings and twin buildings.
\newblock {\em Ann.\ Comb.}, 4:125--137, 2000.

\bibitem{Abramenko/Maldeghem:2001}
Peter Abramenko and Hendrik Van~Maldeghem.
\newblock $1$-twinnings of buildings.
\newblock {\em Math.\ Z.}, 238:187--203, 2001.

\bibitem{Abramenko/Maldeghem:2002}
Peter Abramenko and Hendrik Van~Maldeghem.
\newblock Characterisations of generalized polygons and opposition in rank $2$
  twin buildings.
\newblock {\em J.\ Geom.}, 74:7--28, 2002.

\bibitem{Altmann:2003}
Kristina Altmann.
\newblock The geometry of nondegenerate subspaces of complex orthogonal space.
\newblock Master's thesis, TU Darmstadt, 2003.

\bibitem{Altmann/Gramlich2}
Kristina Altmann and Ralf Gramlich.
\newblock On the hyperbolic unitary geometry.
\newblock Submitted.

\bibitem{Altmann/Gramlich:2006}
Kristina Altmann and Ralf Gramlich.
\newblock On the geometry on the nondegenerate subspaces of nondegenerate
  orthogonal space.
\newblock {\em Bull.\ Belg.\ Math.\ Soc.\ Simon Stevin}, pages 167--179, 2006.

\bibitem{Aschbacher}
Michael Aschbacher.
\newblock Notes on the {Curtis}-{Tits}-{Phan} theorems.
\newblock Unpublished note.

\bibitem{Aschbacher:1977}
Michael Aschbacher.
\newblock A characterization of {Chevalley} groups over fields of odd order,
  parts {I}, {II}.
\newblock {\em Ann.\ of Math.}, 106:353--468, 1977.

\bibitem{Aschbacher:1980}
Michael Aschbacher.
\newblock Correction to: ``{A} characterization of {Chevalley} groups over
  fields of odd order, parts {I}, {II}''.
\newblock {\em Ann.\ of Math.}, 111:411--414, 1980.

\bibitem{Aschbacher:1993}
Michael Aschbacher.
\newblock Simple connectivity of $p$-group complexes.
\newblock {\em Israel J.\ Math.}, 82:1--43, 1993.

\bibitem{Aschbacher:2004}
Michael Aschbacher.
\newblock The status of the classification of the finite simple groups.
\newblock {\em Notices Amer.\ Math.\ Soc.}, 51:736--740, 2004.

\bibitem{Bahls:2005}
Patrick Bahls.
\newblock {\em The isomorphism problem in {Coxeter} groups}.
\newblock Imperial College Press, London, 2005.

\bibitem{Bennett/Gramlich/Hoffman/Shpectorov:2003}
Curtis~D.\ Bennett, Ralf Gramlich, Corneliu Hoffman, and Sergey Shpectorov.
\newblock {Curtis}-{Phan}-{Tits} theory.
\newblock In Alexander~A.\ Ivanov, Martin~W.\ Liebeck, and Jan Saxl, editors,
  {\em Groups, Combinatorics and Geometry: Durham 2001}, pages 13--29, New
  Jersey, 2003. World Scientific.

\bibitem{Bennett/Gramlich/Hoffman/Shpectorov:2007}
Curtis~D.\ Bennett, Ralf Gramlich, Corneliu Hoffman, and Sergey Shpectorov.
\newblock Odd-dimensional orthogonal groups as amalgams of unitary groups, part
  1: general simple connectedness.
\newblock {\em J.\ Algebra}, 312:426--444, 2007.

\bibitem{Bennett/Shpectorov:2004}
Curtis~D.\ Bennett and Sergey Shpectorov.
\newblock A new proof of a theorem of {Phan}.
\newblock {\em J.\ Group Theory}, 7:287--310, 2004.

\bibitem{Benson:1991}
David~J.\ Benson.
\newblock {\em Representations and Cohomology. {II}: Cohomology of Groups and
  Modules}.
\newblock Cambridge University Press, Cambridge, 1991.

\bibitem{Blok/Hoffman:2008}
Rieuwert~J.\ Blok and Corneliu Hoffman.
\newblock A quasi {Curtis}--{Tits}--{Phan} theorem for the symplectic group.
\newblock {\em J.\ Algebra}, 319:4662--4691, 2008.

\bibitem{Blok/Hoffman:2009}
Rieuwert~J.\ Blok and Corneliu Hoffman.
\newblock A {Curtis}--{Tits}--{Phan} theorem for the twin building of type
  $\tilde {A}_{n-1}$.
\newblock {\em J.\ Algebra}, 321:1196--1224, 2009.

\bibitem{Bourbaki:2002}
Nicolas Bourbaki.
\newblock {\em Lie groups and {L}ie algebras. {C}hapters 4--6}.
\newblock Elements of Mathematics (Berlin). Springer-Verlag, Berlin, 2002.
\newblock Translated from the 1968 French original by Andrew Pressley.

\bibitem{Bredon:1993}
Glen~E.\ Bredon.
\newblock {\em Topology and geometry}, volume 139 of {\em Graduate texts in
  mathematics}.
\newblock Springer, Berlin, 1993.

\bibitem{Bridson/Haefliger:1999}
Martin Bridson and Andr{\'e} Haefliger.
\newblock {\em Metric spaces of non-positive curvature}.
\newblock Springer, Berlin, 1999.

\bibitem{Brouwer:1993}
Andries~E.\ Brouwer.
\newblock The complement of a geometric hyperplane in a generalized polygon is
  usually connected.
\newblock In Frank De~Clerck et~al., editor, {\em Finite Geometries and
  Combinatorics}, pages 53--57, 1993.

\bibitem{Brown:1989}
Kenneth~S.\ Brown.
\newblock {\em Buildings}.
\newblock Springer, Berlin, 1989.

\bibitem{Buekenhout:1995}
Francis Buekenhout, editor.
\newblock {\em Handbook of incidence geometry}.
\newblock North-Holland, Amsterdam, 1995.

\bibitem{Buekenhout/Cohen}
Francis Buekenhout and Arjeh~M.\ Cohen.
\newblock Diagram geometry.
\newblock Springer, Berlin, in preparation, {\tt
  http://www.win.tue.nl/{$\sim$}amc/buek}.

\bibitem{Buekenhout/Shult:1974}
Francis Buekenhout and Ernest~E.\ Shult.
\newblock On the foundations of polar geometry.
\newblock {\em Geom.\ Dedicata}, 3:155--170, 1974.

\bibitem{Caprace:2007a}
Pierre-Emmanuel Caprace.
\newblock On $2$-spherical {Kac-Moody} groups and their central extensions.
\newblock {\em Forum\ Math.}, 19(5):763--781, 2007.

\bibitem{Caprace:2009}
Pierre-Emmanuel Caprace.
\newblock {\em ``Abstract'' homomorphisms of split {Kac-Moody} groups}, volume
  198 of {\em Memoirs of the AMS}.
\newblock American Mathematical Society, New York, 2009.

\bibitem{Caprace/Muehlherr:2005}
Pierre-Emmanuel Caprace and Bernhard M{\"u}hlherr.
\newblock Isomorphisms of {Kac-Moody} groups.
\newblock {\em Invent.\ Math.}, 161:361--388, 2005.

\bibitem{Caprace/Muehlherr:2006}
Pierre-Emmanuel Caprace and Bernhard M{\"u}hlherr.
\newblock Isomorphisms of {Kac-Moody} groups which preserve bounded subgroups.
\newblock {\em Adv.\ Math.}, 206:250--278, 2006.

\bibitem{Caprace/Muehlherr:2007}
Pierre-Emmanuel Caprace and Bernhard M{\"u}hlherr.
\newblock Reflection rigidity of $2$-spherical {Coxeter} groups.
\newblock {\em Proc.\ Lond.\ Math.\ Soc.}, 94:520--542, 2007.

\bibitem{Caprace/Remy:2006}
Pierre-Emmanuel Caprace and Bertrand R{\'e}my.
\newblock Simplicit{\'e} abstraite des groupes de {Kac-Moody} non affines.
\newblock {\em C.\ R.\ Acad.\ Sci.\ Paris, Ser.\ I}, 342:539--544, 2006.

\bibitem{Caprace/Remy}
Pierre-Emmanuel Caprace and Bertrand R{\'e}my.
\newblock Simplicity and superrigidity of twin building lattices.
\newblock {\em Invent.\ Math.}, 176:169--221, 2009.

\bibitem{Charney/Davis:2000}
Ruth Charney and Michael Davis.
\newblock When is a {Coxeter} system determined by its {Coxeter} group.
\newblock {\em J.\ Lond.\ Math.\ Soc.}, 61:441--461, 2000.

\bibitem{Chevalley:1966}
Claude Chevalley.
\newblock Certains sch{\'e}mas de groupes semi-simples.
\newblock In {\em 13e S{\'e}minaire N.\ Bourbaki, 1960--1961}. Benjamin, New
  York, 1966.

\bibitem{Cohen:1995}
Arjeh~M.\ Cohen.
\newblock Point-line spaces related to buildings.
\newblock In Buekenhout \cite{Buekenhout:1995}, pages 647--737.

\bibitem{Cohen/Cuypers/Gramlich:2005}
Arjeh~M.\ Cohen, Hans Cuypers, and Ralf Gramlich.
\newblock Local recognition of non-incident point-hyperplane pairs.
\newblock {\em Combinatorica}, 25:271--296, 2005.

\bibitem{Coxeter:1935}
Harold~S.M.\ Coxeter.
\newblock The complete enumeration of finite groups of the form ${R}^2 =
  ({R}_i,{R}_j)^{k_{ij}} = 1$.
\newblock {\em J.\ Lond.\ Math.\ Soc.}, pages 21--25, 1935.

\bibitem{Curtis:1965}
Charles~W.\ Curtis.
\newblock Central extensions of groups of {Lie} type.
\newblock {\em J.\ Reine Angew.\ Math.}, 220:174--185, 1965.

\bibitem{Cuypers:1994}
Hans Cuypers.
\newblock Symplectic geometries, transvection groups, and modules.
\newblock {\em J.\ Comb.\ Theory Ser.\ A}, 65:39--59, 1994.

\bibitem{Das:1994}
Kaustuv~M.\ Das.
\newblock {\em Homotopy and homology of $p$-group complexes}.
\newblock PhD thesis, Caltech, 1994.

\bibitem{Davis:2008}
Michael Davis.
\newblock {\em The geometry and topology of Coxeter groups}.
\newblock London Mathematical Society, 2008.

\bibitem{Delgado/Goldschmidt/Stellmacher:1985}
Alberto Delgado, David~M.\ Goldschmidt, and Bernd Stellmacher.
\newblock {\em Groups and graphs: new results and methods}, volume~6 of {\em
  DMV Seminar}.
\newblock Birkh{\"a}user, Basel, 1985.

\bibitem{Demazure:1965}
Michel Demazure.
\newblock Sch{\'e}mas en groupes r{\'e}ductifs.
\newblock {\em Bull.\ Soc.\ Math.\ France}, 93:369--413, 1965.

\bibitem{Devillers/Gramlich/Muehlherr}
Alice Devillers, Ralf Gramlich, and Bernhard M{\"u}hlherr.
\newblock The sphericity of the complex of nondegenerate subspaces.
\newblock {\em J.\ London Math.\ Soc.}, to appear.

\bibitem{Devillers/Muehlherr:2007}
Alice Devillers and Bernhard M{\"u}hlherr.
\newblock On the simple connectedness of certain subsets of buildings.
\newblock {\em Forum Math.}, 19:955--970, 2007.

\bibitem{Dress/Scharlau:1987}
Andreas Dress and Rudolf Scharlau.
\newblock Gated sets in metric spaces.
\newblock {\em Aequiationes Math.}, 34:112--120, 1987.

\bibitem{Dunlap:2005}
Jonathan~R.\ Dunlap.
\newblock {\em Uniqueness of {Curits}-{Phan}-{Tits} amalgams}.
\newblock PhD thesis, Bowling Green State University, 2005.

\bibitem{Fischer:1971}
Bernd Fischer.
\newblock Finite groups generated by $3$-transpositions, {I}.
\newblock {\em Invent.\ Math.}, 13:232--246, 1971.

\bibitem{Freudenthal:1951}
Hans Freudenthal.
\newblock {\em Oktaven, Ausnahmegruppen und Oktavengeometrie}.
\newblock Mathematisch Instituut der Rijksuniversiteit te Utrecht, Utrecht,
  1951.

\bibitem{Freudenthal:1985}
Hans Freudenthal.
\newblock {Oktaven, Ausnahmegruppen und Oktavengeometrie}.
\newblock {\em Geom.\ Dedicata}, 19:7--63, 1985.

\bibitem{Garland:1973}
Howard Garland.
\newblock $p$-adic curvature and the cohomology of discrete subgroups of
  $p$-adic groups.
\newblock {\em Annals Math.}, 97:375--423, 1973.

\bibitem{Goldschmidt:1980}
David~M.\ Goldschmidt.
\newblock Automorphisms of trivalent graphs.
\newblock {\em Annals Math.}, 111:377--406, 1980.

\bibitem{Gorenstein/Lyons/Solomon:1994}
Daniel Gorenstein, Richard Lyons, and Ronald Solomon.
\newblock {\em The classification of the finite simple groups}, volume 40.1 of
  {\em Mathematical Surveys and Monographs}.
\newblock American Mathematical Society, Providence, 1994.

\bibitem{Gorenstein/Lyons/Solomon:1995}
Daniel Gorenstein, Richard Lyons, and Ronald Solomon.
\newblock {\em The classification of the finite simple groups. Number 2. Part
  I. Chapter G. General group theory}, volume 40.2 of {\em Mathematical Surveys
  and Monographs}.
\newblock American Mathematical Society, Providence, 1995.

\bibitem{Gorenstein/Lyons/Solomon:1998}
Daniel Gorenstein, Richard Lyons, and Ronald Solomon.
\newblock {\em The classification of the finite simple groups. Number 3. Part
  I. Chapter A. Almost simple $K$-groups}, volume 40.3 of {\em Mathematical
  Surveys and Monographs}.
\newblock American Mathematical Society, Providence, 1998.

\bibitem{Gorenstein/Lyons/Solomon:1999}
Daniel Gorenstein, Richard Lyons, and Ronald Solomon.
\newblock {\em The classification of the finite simple groups. Number 4. Part
  II. Chapters 1--4. Uniqueness Theorems}, volume 40.4 of {\em Mathematical
  Surveys and Monographs}.
\newblock American Mathematical Society, Providence, 1999.

\bibitem{Gorenstein/Lyons/Solomon:2002}
Daniel Gorenstein, Richard Lyons, and Ronald Solomon.
\newblock {\em The classification of the finite simple groups. Number 5. Part
  III. Chapters 1--6. The generic case, stages 1--3a}, volume 40.5 of {\em
  Mathematical Surveys and Monographs}.
\newblock American Mathematical Society, Providence, 2002.

\bibitem{Gorenstein/Lyons/Solomon:2006}
Daniel Gorenstein, Richard Lyons, and Ronald Solomon.
\newblock {\em The classification of the finite simple groups. Number 6. Part
  IV. The special odd case}, volume 40.6 of {\em Mathematical Surveys and
  Monographs}.
\newblock American Mathematical Society, Providence, 2006.

\bibitem{Gramlich:2002}
Ralf Gramlich.
\newblock {\em On graphs, geometries, and groups of Lie type}.
\newblock PhD thesis, TU Eindhoven, 2002.
\newblock {\tt http://alexandria.tue.nl/extra2/200211439.pdf}.

\bibitem{Gramlich:2004b}
Ralf Gramlich.
\newblock Line-hyperline pairs of projective spaces and fundamental subgroups
  of linear groups.
\newblock {\em Adv.\ Geom.}, 4:83--103, 2004.

\bibitem{Gramlich:2004c}
Ralf Gramlich.
\newblock On the hyperbolic symplectic geometry.
\newblock {\em J.\ Combin.\ Theory Ser.\ A}, 105:97--110, 2004.

\bibitem{Gramlich:2004d}
Ralf Gramlich.
\newblock {Weak {Phan} systems of type $C_n$}.
\newblock {\em J.\ Algebra}, 280:1--19, 2004.

\bibitem{Gramlich/Hall/Straub}
Ralf Gramlich, Jonathan~I.\ Hall, and Armin Straub.
\newblock The local recognition of reflection graphs of spherical {Coxeter}
  groups.
\newblock submitted.

\bibitem{Gramlich/Hoffman/Nickel/Shpectorov:2005}
Ralf Gramlich, Corneliu Hoffman, Werner Nickel, and Sergey Shpectorov.
\newblock Even-dimensional orthogonal groups as amalgams of unitary groups.
\newblock {\em J.\ Algebra}, pages 141--173, 2005.

\bibitem{Gramlich/Hoffman/Shpectorov:2003}
Ralf Gramlich, Corneliu Hoffman, and Sergey Shpectorov.
\newblock A {Phan}-type theorem for {${\rm Sp}(2n,q)$}.
\newblock {\em J.\ Algebra}, 264:358--384, 2003.

\bibitem{Gramlich/Horn/Muehlherr}
Ralf Gramlich, Max Horn, and Bernhard M{\"u}hlherr.
\newblock Abstract involutions of algebraic groups and of {Kac}--{Moody}
  groups.
\newblock submitted.

\bibitem{Gramlich/Horn/Nickel:2006}
Ralf Gramlich, Max Horn, and Werner Nickel.
\newblock The complete {Phan}-type theorem for {${\rm Sp}(2n,q)$}.
\newblock {\em J.\ Group Theory}, 9:603--626, 2006.

\bibitem{Gramlich/Horn/Nickel:2007}
Ralf Gramlich, Max Horn, and Werner Nickel.
\newblock Odd-dimensional orthogonal groups as amalgams of unitary groups, part
  2: machine computations.
\newblock {\em J.\ Algebra}, 316:591--607, 2007.

\bibitem{Gramlich/Horn/Pasini/Maldeghem}
Ralf Gramlich, Max Horn, Antonio Pasini, and Hendrik Van~Maldeghem.
\newblock Intransitive geometries and fused amalgams.
\newblock {\em J.\ Group Theory}, 11:443--464, 2008.

\bibitem{Gramlich/Muehlherr:2008}
Ralf Gramlich and Bernhard M{\"u}hlherr.
\newblock Lattices from involutions of {Kac-Moody} groups.
\newblock {\em Oberwolfach reports}, 5:139--140, 2008.

\bibitem{Gramlich/Maldeghem:2006}
Ralf Gramlich and Hendrik Van~Maldeghem.
\newblock Intransitive geometries.
\newblock {\em Proc.\ Lond.\ Math.\ Soc.}, 93:666--692, 2006.

\bibitem{Gramlich/Witzel}
Ralf Gramlich and Stefan Witzel.
\newblock The sphericity of generalized {Phan} geometries of type {$B_n$} and
  {$C_n$} and the {Phan}-type theorem of type {$F_4$}.
\newblock submitted.

\bibitem{Gross/Tucker:1987}
Jonathan~L.\ Gross and Thomas~W.\ Tucker.
\newblock {\em Topological graph theory}.
\newblock Wiley, New York, 1987.

\bibitem{Grundhoefer:2002}
Theo Grundh{\"o}fer.
\newblock Basics on buildings.
\newblock In Tent \cite{Tent:2002}, pages 1--21.

\bibitem{Hall:1988}
Jonathan~I.\ Hall.
\newblock The hyperbolic lines of finite symplectic spaces.
\newblock {\em J.\ Combin.\ Theory Ser.\ A}, 47:284--298, 1988.

\bibitem{Hall:1989}
Jonathan~I.\ Hall.
\newblock Graphs, geometry, $3$-transpositions, and symplectic
  {$\mathbb{F}_2$}-transvection groups.
\newblock {\em Proc.\ London Math.\ Soc.}, 58:89--111, 1989.

\bibitem{Horn:2005}
Max Horn.
\newblock Amalgams of unitary groups in {$Sp(2n,q)$}.
\newblock Master's thesis, Tech\-nische Universit{\"a}t Darmstadt, 2005.

\bibitem{Horn:2008b}
Max Horn.
\newblock {\em Involutions of Kac-Moody groups}.
\newblock PhD thesis, TU Darmstadt, 2008.

\bibitem{Horn:2008a}
Max Horn.
\newblock On the {Phan} system of the {Schur} cover of {$\mathrm{SU}(4,3^2)$}.
\newblock {\em Des.\ Codes Cryptogr.}, 47(1-3):243--247, 2008.

\bibitem{Hughes/Piper:1973}
Daniel~R.\ Hughes and Fred~C.\ Piper.
\newblock {\em Projective Planes}.
\newblock Springer, Berlin, 1973.

\bibitem{Humphreys:1972}
James~E.\ Humphreys.
\newblock Remarks on ``{A} theorem on special linear groups''.
\newblock {\em J.\ Algebra}, 22:316--318, 1972.

\bibitem{Humphreys:1987}
James~E.\ Humphreys.
\newblock The {Steinberg} representation.
\newblock {\em Bull.\ Amer.\ Math.\ Soc.\ (N.S.)}, 16:247--263, 1987.

\bibitem{Humphreys:1990}
James~E.\ Humphreys.
\newblock {\em Reflection groups and {Coxeter} groups}.
\newblock Cambridge University Press, Cambridge, 1990.

\bibitem{Ivanov:1999}
Alexander~A.\ Ivanov.
\newblock {\em Geometry of sporadic groups {I}: Petersen and tilde geometries},
  volume~76 of {\em Encyclopedia of Mathematics}.
\newblock Cambridge University Press, 1999.

\bibitem{Ivanov/Shpectorov:2002}
Alexander~A.\ Ivanov and Sergey Shpectorov.
\newblock {\em Geometry of sporadic groups {II}: representations and amalgams},
  volume~91 of {\em Encyclopedia of Mathematics}.
\newblock Cambridge University Press, 2002.

\bibitem{Klein:1872}
Felix Klein.
\newblock {\em Vergleichende Betrachtungen \"uber neuere geometrische
  For\-schungen}.
\newblock Verlag von Andreas Deichert, Erlangen, 1872.
\newblock Also known as {\em Erlanger Programm}.

\bibitem{Kuesefoglu:1979}
Ayse~S.\ K{\"u}sefoglu.
\newblock The second degree cohomology of finite orthogonal groups. {I}.
\newblock {\em J.\ Algebra}, 56:207--220, 1979.

\bibitem{Kuesefoglu:1980}
Ayse~S.\ K{\"u}sefoglu.
\newblock The second degree cohomology of finite orthogonal groups. {II}.
\newblock {\em J.\ Algebra}, 67:88--109, 1980.

\bibitem{Margulis:1991}
Gregori Margulis.
\newblock {\em Discrete subgroups of semisimple {L}ie groups}.
\newblock Springer, Berlin, 1991.

\bibitem{Muehlherr:1994}
Bernhard M{\"u}hlherr.
\newblock {\em Some Contributions to the Theory of Buildings Based on the Gate
  Property}.
\newblock PhD thesis, Uni T{\"u}bingen, 1994.

\bibitem{Muehlherr}
Bernhard M{\"u}hlherr.
\newblock On the simple connectedness of a chamber system associated to a twin
  building.
\newblock Unpublished note, 1996.

\bibitem{Muehlherr:1998}
Bernhard M{\"u}hlherr.
\newblock A rank $2$ characterization of twinnings.
\newblock {\em European J.\ Combin.}, 19:603--612, 1998.

\bibitem{Muehlherr:1999}
Bernhard M{\"u}hlherr.
\newblock Locally split and locally finite twin buildings of $2$-spherical
  type.
\newblock {\em J.\ Reine Angew.\ Math.}, 511:119--143, 1999.

\bibitem{Muehlherr:2002}
Bernhard M{\"u}hlherr.
\newblock Twin buildings.
\newblock In Tent \cite{Tent:2002}, pages 103--117.

\bibitem{Pasini:1985}
Antonio Pasini.
\newblock Some remarks on covers and apartments.
\newblock In Catharine~A.\ Baker and Lynn~M.\ Batten, editors, {\em Finite
  geometries}, volume 103 of {\em Lecture Notes in Pure and Applied
  Mathematics}, pages 223--250. Dekker, New York, 1985.

\bibitem{Pasini:1994}
Antonio Pasini.
\newblock {\em Diagram Geometries}.
\newblock Clarendon Press, Oxford, 1994.

\bibitem{Payne/Thas:1984}
Stanley~E.\ Payne and Joseph~A.\ Thas.
\newblock {\em Finite generalized quadrangles}.
\newblock Pitman, Boston, 1984.

\bibitem{Phan:1970}
Kok{-}Wee Phan.
\newblock A theorem on special linear groups.
\newblock {\em J.\ Algebra}, 16:509--518, 1970.

\bibitem{Phan:1977}
Kok{-}Wee Phan.
\newblock On groups genererated by three-dimensional special unitary groups,
  {I}.
\newblock {\em J.\ Austral.\ Math.\ Soc.\ Ser.\ A}, 23:67--77, 1977.

\bibitem{Phan:1977a}
Kok{-}Wee Phan.
\newblock On groups genererated by three-dimensional special unitary groups,
  {II}.
\newblock {\em J.\ Austral.\ Math.\ Soc.\ Ser.\ A}, 23:129--146, 1977.

\bibitem{Pickert:1975}
G{\"u}nter Pickert.
\newblock {\em Projektive Ebenen}.
\newblock Springer, Berlin, second edition, 1975.

\bibitem{Remy:2002}
Bertrand R{\'e}my.
\newblock {\em Groupes de {Kac-Moody} d\'eploy\'es et presque d\'eploy\'es},
  volume 277 of {\em Ast\'erisque}.
\newblock Soci{\'e}t{\'e} Math{\'e}matiques de France, Paris, 2002.

\bibitem{Roberts:2005}
Adam~E.\ Roberts.
\newblock {\em A {Phan}-type theorem for orthogonal groups}.
\newblock PhD thesis, Bowling Green State University, 2005.

\bibitem{Ronan:1989}
Mark~A.\ Ronan.
\newblock {\em Lectures on buildings}, volume~7 of {\em Perspectives in
  Mathematics}.
\newblock Academic Press, Boston, 1989.

\bibitem{Ronan:2000}
Mark~A.\ Ronan.
\newblock Local isometries of twin buildings.
\newblock {\em Math.\ Z.}, 234:435--455, 2000.

\bibitem{Ronan:2002}
Mark~A.\ Ronan.
\newblock Twin trees and twin buildings.
\newblock In Tent \cite{Tent:2002}, pages 119--137.

\bibitem{Ronan/Tits:1994}
Mark~A.\ Ronan and Jacques Tits.
\newblock Twin trees. {I}.
\newblock {\em Invent.\ Math.}, 116:463--479, 1994.

\bibitem{Scharlau:1995}
Rudolf Scharlau.
\newblock Buildings.
\newblock In Buekenhout \cite{Buekenhout:1995}, pages 477--645.

\bibitem{Seifert/Threlfall:1934}
Herbert Seifert and William Threlfall.
\newblock {\em Lehrbuch der Topologie}.
\newblock Chelsea Publishing Company, New York, 1934.

\bibitem{Serre:1977}
Jean-Pierre Serre.
\newblock {\em Arbres, amalgames, $SL_2$}, volume~46 of {\em Ast\'erisque}.
\newblock Soc.\ Math.\ France, Paris, 1977.

\bibitem{Serre:2003}
Jean-Pierre Serre.
\newblock {\em Trees}.
\newblock Springer, Berlin, 2003.
\newblock Corrected second printing.

\bibitem{Solomon:1969}
Louis Solomon.
\newblock The {Steinberg} character of a finite group with {$BN$}-pair.
\newblock In {\em Theory of finite groups}, pages 213--221. Benjamin, New York,
  1969.

\bibitem{Spanier:1966}
Edwin~H.\ Spanier.
\newblock {\em Algebraic Topology}.
\newblock McGraw-Hill, New York, 1966.

\bibitem{Steinberg:1962}
Robert Steinberg.
\newblock G{\'e}n{\'e}rateurs, relations, et rev{\^e}tements de groupes
  al\-g{\'e}\-briques.
\newblock In {\em Colloque sur la th{\'e}orie des groupes alg{\'e}briques},
  pages 113--127. Bruxelles, 1962.

\bibitem{Steinberg:1968}
Robert Steinberg.
\newblock Lectures on {Chevalley} groups.
\newblock Mimeographed lecture notes, Yale University, New Haven, 1968.

\bibitem{Stroppel:1992}
Markus Stroppel.
\newblock Reconstruction of incidence geometries from groups of automorphisms.
\newblock {\em Arch.\ Math.}, 58:621--624, 1992.

\bibitem{Stroppel:1993}
Markus Stroppel.
\newblock A categorial glimpse at the reconstruction of geometries.
\newblock {\em Geom.\ Dedicata}, 46:47--60, 1993.

\bibitem{Stroth:1989}
Gernot Stroth.
\newblock Some geometries for ${M}c{L}$.
\newblock {\em Comm.\ Alg.}, 17:2825--2833, 1989.

\bibitem{Tent:2002}
Katrin Tent, editor.
\newblock {\em $BN$-pairs and groups of finite Morley rank. Tits buildings and
  the model theory of groups (W{\"u}rzburg, 2000)}.
\newblock Cambridge University Press, Cambridge, 2002.

\bibitem{Timmesfeld:1991}
Franz~Georg Timmesfeld.
\newblock Groups generated by $k$-root subgroups.
\newblock {\em Invent.\ Math.}, 106:575--666, 1991.

\bibitem{Timmesfeld:1998}
Franz~Georg Timmesfeld.
\newblock Presentations for certain {Chevalley} groups.
\newblock {\em Geom.\ Dedicata}, 73:85--117, 1998.

\bibitem{Timmesfeld:1999}
Franz~Georg Timmesfeld.
\newblock Abstract root subgroups and quadratic action.
\newblock {\em Adv.\ Math.}, 142:1--150, 1999.

\bibitem{Timmesfeld:2001}
Franz~Georg Timmesfeld.
\newblock {\em Abstract root subgroups and simple groups of Lie type},
  volume~95 of {\em Monographs in Mathematics}.
\newblock Birkh\"auser, Basel, 2001.

\bibitem{Timmesfeld:2003}
Franz~Georg Timmesfeld.
\newblock On the {Steinberg}-presentation for {Lie}-type groups.
\newblock {\em Forum Math.}, 15:645--663, 2003.

\bibitem{Timmesfeld:2004}
Franz~Georg Timmesfeld.
\newblock The {Curtis}-{Tits} presentation.
\newblock {\em Adv.\ Math.}, 189:38--67, 2004.

\bibitem{Timmesfeld:2006}
Franz~Georg Timmesfeld.
\newblock Steinberg-type presentations for {Lie}-type groups.
\newblock {\em J.\ Algebra}, 300:806--819, 2006.

\bibitem{Tits:1981}
Jacques Tits.
\newblock A local approach to buildings.
\newblock In {\em The Geometric Vein -- The Coxeter Festschrift}, pages
  519--547, New York. Springer.

\bibitem{Tits:1962}
Jacques Tits.
\newblock Groupes semi-simples isotropes.
\newblock In {\em Colloque sur la th{\'e}orie des groupes alg{\'e}briques},
  pages 137--146, Bruxelles, 1962.

\bibitem{Tits:1974}
Jacques Tits.
\newblock {\em Buildings of spherical type and finite {$BN$}-pairs}, volume 386
  of {\em Lecture Notes in Mathematics}.
\newblock Springer, Berlin, 1974.

\bibitem{Tits:1986}
Jacques Tits.
\newblock Ordonn\'es, immeubles et sommes amalgam\'ees.
\newblock {\em Bull.\ Soc.\ Math.\ Belg.\ S\'er.\ A}, 38:367--387, 1986.

\bibitem{Tits:1987}
Jacques Tits.
\newblock Uniqueness and presentation of {Kac-Moody} groups over fields.
\newblock {\em J.\ Algebra}, 105:542--573, 1987.

\bibitem{Tits:1989}
Jacques Tits.
\newblock R{\'e}sum{\'e} de cours.
\newblock In {\em Annuaire du Coll{\`e}ge de France, 89e ann{\'e}e, 1989 --
  1990}, pages 87--103. Coll{\`e}ge de France, 1990.

\bibitem{Tits:1992}
Jacques Tits.
\newblock {Twin buildings and groups of Kac-Moody type}.
\newblock In Martin~W.\ Liebeck and Jan Saxl, editors, {\em Groups,
  Combinatorics and Geometry}, volume 165 of {\em LMS Lecture Note Series},
  pages 249--286, Cambridge, 1992. Cambridge University Press.

\bibitem{Tits/Weiss:2002}
Jacques Tits and Richard Weiss.
\newblock {\em Moufang Polygons}.
\newblock Springer, Berlin, 2002.

\bibitem{tomDieck:2000}
Tammo tom Dieck.
\newblock {\em Topologie}.
\newblock Walter de Gruyter, Berlin, second edition, 2000.

\bibitem{Maldeghem:1998}
Hendrik Van~Maldeghem.
\newblock {\em Generalized Polygons}, volume~93 of {\em Monographs in
  Mathematics}.
\newblock Birkh{\"a}user, Basel, 1998.

\bibitem{Veldkamp:1959}
Ferdinand~D.\ Veldkamp.
\newblock Polar geometry {I-IV}.
\newblock {\em Proc.\ Nederl.\ Akad.\ Wet.\ Ser.\ A}, 62:515--551, 1959.

\bibitem{Veldkamp:1960}
Ferdinand~D.\ Veldkamp.
\newblock Polar geometry {V}.
\newblock {\em Proc.\ Nederl.\ Akad.\ Wet.\ Ser.\ A}, 63:207--212, 1960.

\bibitem{Weiss:2003}
Richard Weiss.
\newblock {\em The structure of spherical buildings}.
\newblock Princeton University Press, Princeton, 2003.

\bibitem{Wich:1996}
Anke Wich.
\newblock {Skizzierte Geometrien}.
\newblock Master's thesis, TU Darmstadt, 1996.

\bibitem{Wong:1974}
Warren~J.\ Wong.
\newblock Generators and relations for classical groups.
\newblock {\em J.\ Algebra}, 32:529--553, 1974.

\end{thebibliography}

}

\end{document}